\documentclass[a4paper,fleqn]{cas-sc}

\usepackage[authoryear,longnamesfirst]{natbib}

\def\tsc#1{\csdef{#1}{\textsc{\lowercase{#1}}\xspace}}
\tsc{WGM}
\tsc{QE}

\usepackage{amssymb}
\usepackage{latexsym}
\usepackage{amsmath}
\usepackage{subcaption}
\usepackage{float}
\usepackage{siunitx}
\usepackage{placeins}

\DeclareMathOperator{\dive}{\nabla\cdot}

\begin{document}
\let\WriteBookmarks\relax
\def\floatpagepagefraction{1}
\def\textpagefraction{.001}

\shorttitle{Monotonization approach for DG}    

\shortauthors{G. Orlando}  

\title[mode = title]{A filtering monotonization approach for DG discretizations of hyperbolic problems}  

%
\author[1]{Giuseppe Orlando}[orcid = 0000-0002-7119-4231]

\cormark[1]

\ead{giuseppe.orlando@polimi.it}

\affiliation[1]{organization={MOX, Dipartimento di Matematica, Politecnico di Milano},
                addressline={Piazza Leonardo da Vinci 32}, 
                city={Milano},
                postcode={20133}, 
                state={},
                country={Italy}}

\cortext[1]{Corresponding author}

\nonumnote{}

\begin{abstract}
We introduce a filtering technique for Discontinuous Galerkin approximations of hyperbolic problems. Following an approach already proposed for the Hamilton-Jacobi equations by other authors, we aim at reducing the spurious oscillations that arise in presence of discontinuities when high order spatial discretizations are employed. This goal is achieved using a filter function that keeps the high order scheme when the solution is regular and switches to a monotone low order approximation if it is not. The method has been implemented in the framework of the \textit{deal.II} numerical library, whose mesh adaptation capabilities are also used to reduce the region in which the low order approximation is used. A number of numerical experiments demonstrate the potential of the proposed filtering technique.
\end{abstract}




\maketitle

\section{Introduction}
\label{sec:intro} \indent

The Discontinuous Galerkin (DG) method has proven itself a very valuable tool for applications to computational fluid dynamics problems in a great variety of flow regimes, see e.g. the seminal contributions \cite{bassi:1997a, bassi:1997b, cockburn:1989a, cockburn:1989b, cockburn:1990, cockburn:1991, cockburn:1998} as well as the reviews in \cite{giraldo:2020, karniadakis:2005}, among many others. For hyperbolic problems, however, spurious oscillations can arise around shocks and other discontinuities when high order spatial discretizations are used. Furthermore, in many applications, maintaining non negativity of the numerical solutions is essential to preserve their correct physical meaning. In order to address these well known issues, a number of monotonization techniques have been proposed in the literature for DG methods. While a full survey of this topic goes beyond the scope of the paper, we review briefly here some of the most popular techniques. 

In general, monotonization techniques for DG methods have been inherited from finite difference and finite volume approaches. For example, starting with \cite{cockburn:1989a, cockburn:1989b, cockburn:1991}, slope limiting techniques have been employed, while other authors have investigated WENO methods \cite{shu:2003, shu:2016} and flux corrected transport methods \cite{kuzmin:2002, restelli:2006, kuzmin:2012}. Another approach is based on the identification of the regions where the discontinuities are located, in which the mesh is then refined and/or the order of the spatial discretization is lowered in order to exploit the monotonicity of (most) low order approximations. In recent years, the very successful MOOD approach has been proposed in \cite{dumbser:2016, dumbser:2014, loubere:2014, zanotti:2015}, which is also based on the identification of the regions of discontinuity and on the switch from a high order DG method to a monotonic first order finite volume method on a locally refined mesh built around the quadrature nodes used by the DG method.

The method proposed in this paper is inspired by the filtering approach outlined in \cite{bokanowski:2016, sahu:2015}. More specifically, a filter function is employed in such a way that, where the solution is regular, we keep the high order solution, whereas otherwise we switch to a low order method. While the proposed strategy is conceptually similar to that of the MOOD approach, the main novelty of the proposed method is that we do not rely on a regularity indicator and that a monotonic solution is retrieved (almost) automatically. For the spatial discretization, we use the DG approach implemented in the numerical library \textit{deal.II} \cite{bangerth:2007}, which provides $h-$refinement capabilities that are exploited in order to reduce the size of the region where the low order approximation is applied. 

The model problem is introduced in Section \ref{sec:model}, along with the space and time discretizations that will be employed. The filtering monotonization approach is introduced in Section \ref{sec:outline} for a scalar hyperbolic PDE and extended in Section \ref{sec:euler} to the inviscid Euler equations. Numerical results validating the proposed approach are presented in Section \ref{sec:numerical}, while some conclusions and perspectives for future work are presented in Section \ref{sec:conclu}.

\section{Model problem and discretization}
\label{sec:model} 

We consider as model problem the nonlinear conservation law
\begin{equation}
\label{eq:model}
\frac{\partial u}{\partial t} + \dive \mathbf{F}(u) = 0,
\end{equation}
where $ \mathbf{F}(u) $ denotes a $d-$dimensional vector field that depends on the unknown $ u = u(\mathbf{x},t), $ and $ \mathbf{x}\in\mathbb{R}^d, $ generally in a non linear way. Simple examples are the linear advection equation and the Burgers equation. We consider a decomposition of the domain \(\Omega\) into a family of quadrilaterals \(\mathcal{T}_h\), where each element is denoted by \(K\). The skeleton \(\mathcal{E}\) denotes the set of all element faces and \(\mathcal{E} = \mathcal{E}^{I} \cup \mathcal{E}^{B}\), where \(\mathcal{E}^{I}\) is the subset of interior faces and \(\mathcal{E}^{B}\) is the subset of boundary faces. Suitable jump and average operators can then be defined as customary for Discontinuous Galerkin discretizations. A face \(\Gamma \in \mathcal{E}^{I}\) shares two elements that we denote by \(K^{+}\) with outward unit normal \(\mathbf{n}^{+}\) and \(K^{-}\) with outward unit normal \(\mathbf{n}^{-}\), whereas for a face \(\Gamma \in \mathcal{E}^{B}\) we denote by \(\mathbf{n}\) the outward unit normal. For a scalar function \(\varphi\)  the jump is defined as
\[\left[\left[\varphi\right]\right] = \varphi^{+}\mathbf{n}^{+} + \varphi^{-}\mathbf{n}^{-} \quad \text{if }\Gamma \in \mathcal{E}^{I} \qquad \left[\left[\varphi\right]\right] = \varphi\mathbf{n} \quad \text{if }\Gamma \in \mathcal{E}^{B}.\]
The average is defined as
\[\left\{\left\{\varphi\right\}\right\} = \frac{1}{2}\left(\varphi^{+} + \varphi^{-}\right) \quad \text{if }\Gamma \in \mathcal{E}^{I} \qquad \left\{\left\{\varphi\right\}\right\} = \varphi \quad \text{if }\Gamma \in \mathcal{E}^{B}.\]
Similar definitions apply for a vector function \(\boldsymbol{\varphi}\):
\begin{align*}
&\left[\left[\boldsymbol{\varphi}\right]\right] = \boldsymbol{\varphi}^{+}\cdot\mathbf{n}^{+} 
+\boldsymbol{\varphi}^{-}\cdot\mathbf{n}^{-} \quad \text{if }\Gamma \in \mathcal{E}^{I} \qquad 
\left[\left[\boldsymbol{\varphi}\right]\right] = \boldsymbol{\varphi}\cdot\mathbf{n} \quad \text{if }\Gamma \in \mathcal{E}^{B} \\
&\left\{\left\{\boldsymbol{\varphi}\right\}\right\} = \frac{1}{2}\left(\boldsymbol{\varphi}^{+} + \boldsymbol{\varphi}^{-}\right) \quad \text{if }\Gamma \in \mathcal{E}^{I} \qquad \left\{\left\{\boldsymbol{\varphi}\right\}\right\} = \boldsymbol{\varphi} \quad \text{if }\Gamma \in \mathcal{E}^{B}.
\end{align*}
We also introduce the following finite element spaces
\begin{equation}\label{eq:fe_spaces}
Q_k = \left\{v \in L^2(\Omega) : v\rvert_K \in \mathbb{Q}_k \quad \forall K \in \mathcal{T}_h\right\},
\end{equation}
where \(\mathbb{Q}_k\) is the space of polynomials of degree \(k\) in each coordinate direction. The spatial discretization coincides with that described in \cite{arndt:2022} and implemented in the \textit{deal.II} library, so that it does not introduce any particular novelty. The shape functions correspond to the products of Lagrange interpolation polynomials for the support points of \(\left(k + 1\right)\)-order Gauss-Lobatto quadrature rule in each coordinate direction.

Concerning the time discretization, we will consider here only the well known TVD Runge Kutta methods described in \cite{gottlieb:1998, gottlieb:2001}. These are high order time discretization schemes that preserve the strong stability properties of first order explicit Euler time stepping and are known as Strong Stability Preserving (SSP) methods. For the convenience of the reader, we briefly recall here the second order and the third order optimal SSP Runge-Kutta methods derived in \cite{gottlieb:1998} for an ordinary differential equation \(y^{\prime} = \mathcal{N}(y)\). The second order scheme reads as follows:
\begin{align}
v^{(1)} &= v^{n} + \Delta t\mathcal{N}(v^{n}) \\
v^{n+1} &= \frac{1}{2} v^{n} + \frac{1}{2} v^{(1)} + \frac{1}{2}\Delta t\mathcal{N}(v^{(1)}),
\end{align}
where \(v^{n} \approx \mathit{y}(t^{n})\) and \(\Delta t\) denotes the time step. The third order method is given instead by:
\begin{align}
v^{(1)} &= v^{n} + \Delta t\mathcal{N}(v^{n}) \\
v^{(2)} &= \frac{3}{4} v^{n} + \frac{1}{4}v^{(1)} + \frac{1}{4}\Delta t\mathcal{N}(v^{(1)}) \\
v^{n+1} &= \frac{1}{3} v^{n} + \frac{2}{3}v^{(2)} + \frac{2}{3}\Delta t\mathcal{N}(v^{(2)}).
\end{align}
Each stage of the TVD method can be represented as 
\begin{equation}
\mathbf{u} = \mathbf{S}(\mathbf{v}),
\end{equation}
where $ \mathbf{u},\mathbf{v} $ denote the new and old values, respectively, of the vector containing the discrete degrees of freedom which identify the spatial approximation to the solutions of \eqref{eq:model}. $ \mathbf{S} $ denotes formally the solution operator associated to a specific time and space discretization. The transition from  $ \mathbf{v} $ to $ \mathbf{u} $ can be interpreted as an advancement in time of $ \alpha \Delta t $ time units, where $ \alpha $ depends on the details of the TVD method and on the specific stage considered. We will denote by $ \mathbf{S}^M $ the discrete operator associated to the monotonic, low order spatial discretization and by $ \mathbf{S}^H $ that associated to a high order, not monotonic spatial discretization.

\section{Outline of filtering monotonization approach}
\label{sec:outline}

We will now introduce the application of the filtering approach proposed in
\cite{bokanowski:2016, sahu:2015} in the above outlined context. First of all, a filter function \(F\) must be introduced. This can be defined in several ways, for example
\begin{equation} \label{eq:filter_function_a}
F_1(x) = x\mathbf{1}_{\left|x\right|\le 1},
\end{equation}
which corresponds to the Oberman-Salvador filter function originally employed in \cite{oberman:2015} or
\begin{equation} \label{eq:filter_function_b}
F_2(x) = \text{sign}(x)\max\left(1 - \left|\left|x\right| - 1\right|,0\right),
\end{equation}
the so-called Froese and Oberman's filter function originally introduced in \cite{froese:2013}. In the simplest possible filtering approach, the filtered version of $ \mathbf{u} $ can be defined as
\begin{equation} \label{eq:filter}
\mathbf{u}^F = \mathbf{S}^M(\mathbf{v}) + \varepsilon\alpha\Delta t
F\left(\frac{\mathbf{S}^H(\mathbf{v}) - \mathbf{S}^M(\mathbf{v})}{\varepsilon\alpha\Delta t}\right),
\end{equation}
where the low order solution \(\mathbf{S}^{M}\) is computed on the nodes of the high order solution \(\mathbf{S}^H\) and \(\varepsilon\) is a suitable parameter, depending on the time and space discretization parameters, such that
\[\lim_{\left(\Delta t, h\right)\to 0}\epsilon(\Delta t,h) = 0,\] 
with \(h = \max\{\mathrm{diam}(K) | K \in \mathcal{T}_h \}\). More details about the choice of \(\varepsilon\) will be given in Section \ref{sec:numerical}. Notice that the filter function is applied componentwise. In this way, as discussed in \cite{bokanowski:2016}, the high order method is only applied to the components $ i $ for which
$$
\frac{|\mathbf{S}^H(\mathbf{v})_i-\mathbf{S}^M(\mathbf{v})_i|}{\varepsilon\alpha\Delta t} \le 1, \quad i = 1,...,\text{dim}(Q_{h}).
$$
As explained in \cite{bokanowski:2016}, \(\epsilon\) has to be chosen in such a way that 
\[\varepsilon \ge c_0 h,\]
where \(c_0\) is a sufficiently large constant. As we will see in Section \ref{sec:numerical}, the aforementioned approach is very dissipative and, unless a very large value of \(c_0\) is adopted, it yields solutions that essentially coincide with the low order one. Therefore, we propose the alternative filtering strategy

\begin{align} \label{eq:filter_rel}
\mathbf{u}^F_{i} &= \mathbf{S}^M(\mathbf{v})_{i} \nonumber \\
&+ \beta \mathbf{S}^M(\mathbf{v})_i
F\left(\frac{\mathbf{S}^H(\mathbf{v})_{i} - \mathbf{S}^M(\mathbf{v})_{i}}{\beta \mathbf{S}^M(\mathbf{v})_i}\right), \quad i = 1,...,\text{dim}(Q_k),
\end{align}
where \(\beta > 0\) is a suitable parameter that represents a tolerance for the ``componentwise relative difference'' \(\frac{\mathbf{S}^H(\mathbf{v})_{i} - \mathbf{S}^M(\mathbf{v})_{i}}{\mathbf{S}^M(\mathbf{v})_i}\), so that when \(\left|\frac{\mathbf{S}^H(\mathbf{v})_{i} - \mathbf{S}^M(\mathbf{v})_{i}}{\mathbf{S}^M(\mathbf{v})_i}\right| \le \beta\), we resort to the high order solution. Also in this case, a too small value of \(\beta\) provides results that are in practice coincident with the low order solution. Appropriate choices for \(\beta\) will also be presented in Section \ref{sec:numerical}. 

\section{Extension to the Euler equations}
\label{sec:euler} 

In this section we present the extension of the strategy \eqref{eq:filter_rel} to the Euler equations 
\begin{equation}
\frac{\partial \mathbf{w}}{\partial t} + \dive\mathbf{F}\left(\mathbf{w}\right) = \mathbf{0}, 
\end{equation} 
where
\[\mathbf{w} = \begin{bmatrix}
\rho \\
\rho \mathbf{u} \\
\rho E
\end{bmatrix} \qquad 
\mathbf{F} = \begin{bmatrix}
\rho \mathbf{u} \\
\rho \mathbf{u} \otimes \mathbf{u} + p\mathbf{I} \\
\left(\rho E + p\right)\mathbf{u}
\end{bmatrix}.\]
Here, \(\rho\) is the fluid density, \(\mathbf{u}\) is the fluid velocity, \(p\) is the pressure and \(E\) is the total energy per unit of mass and \(\mathbf{I}\) is the \(d\)-dimensional identity matrix. The above equations must be complemented by an equation of state (EOS). In this work we consider the classical ideal gas EOS
\begin{equation}
p = \left(\gamma - 1\right)\left(\rho E - \frac{1}{2}\rho \mathbf{u} \cdot \mathbf{u}\right),
\end{equation}
where \(\gamma = 1.4\) is the specific heats ratio. Hence, the application of the filtering approach to the density reads as follows:
\begin{align}\label{eq:filter_rel_rho}
\boldsymbol{\rho}^{F}_{i} &= \mathbf{S}^{M}(\boldsymbol{\rho})_{i} \nonumber \\
&+ \beta_{\rho} \mathbf{S}^{M}(\boldsymbol{\rho})_{i} F\left(\frac{\mathbf{S}^{H}(\boldsymbol{\rho})_{i} - \mathbf{S}^{M}(\boldsymbol{\rho})_{i}}{\beta_{\rho} \mathbf{S}^{M}(\boldsymbol{\rho})_{i}}\right), \quad i = 1,...,\text{dim}(Q_k),
\end{align} 
where \(\boldsymbol{\rho}\) is the vector of the degrees of freedom for the density and \(\beta_{\rho}\) is the tolerance parameter for the density. In analogy to \cite{loubere:2014}, we choose to perform the filtering procedure for all the conserved variables, namely also for \(\rho \mathbf{u}\) and \(\rho E\), and their formulation is analogous to \eqref{eq:filter_rel_rho}:
\begin{align}
\boldsymbol{\rho}\mathbf{u}^{F}_{i} &= \mathbf{S}^{M}(\boldsymbol{\rho}\mathbf{u})_{i} \nonumber \\
&+ \beta_{\rho\mathbf{u}} \mathbf{S}^{M}(\boldsymbol{\rho}\mathbf{u})_{i} F\left(\frac{\mathbf{S}^{H}(\boldsymbol{\rho}\mathbf{u})_{i} - \mathbf{S}^{M}(\boldsymbol{\rho}\mathbf{u})_{i}}{\beta_{\rho\mathbf{u}} \mathbf{S}^{M}(\boldsymbol{\rho}\mathbf{u})_{i}}\right), \quad i = 1,...,\left[\text{dim}(Q_k)\right]^{d} \\
\boldsymbol{\rho}\mathbf{E}^{F}_{i} &= \mathbf{S}^{M}(\boldsymbol{\rho}\mathbf{E})_{i} \nonumber \\
&+ \beta_{\rho E} \mathbf{S}^{M}(\boldsymbol{\rho}\mathbf{E})_{i} F\left(\frac{\mathbf{S}^{H}(\boldsymbol{\rho}\mathbf{E})_{i} - \mathbf{S}^{M}(\boldsymbol{\rho}\mathbf{E})_{i}}{\beta_{\rho E} \mathbf{S}^{M}(\boldsymbol{\rho}\mathbf{E})_{i}}\right), \quad i = 1,...,\text{dim}(Q_k).
\end{align}

\section{Numerical experiments}
\label{sec:numerical} 

The numerical scheme outlined in the previous Sections has been validated in a number of  benchmarks. We set 
\[{\cal H} = \min\{\mathrm{diam}(K) | K \in \mathcal{T}_h\}\]
and we define the Courant number:
\begin{equation}
\label{eq:Courant}
C = kU \Delta t/{\cal H},
\end{equation}
where \(U\) is the magnitude of the flow velocity. In the case of the Euler equations, the Courant number \(C\) is defined as:
\begin{equation}
\label{eq:Courant_Euler}
C = k\left(U + c\right)\Delta t/{\cal H},
\end{equation}
where \(c = \sqrt{\gamma \frac{p}{\rho}}\) is the speed of sound. We chose to employ mainly \(k = 1\) and \(k = 2\) in combination with the second order SSP and the third order SSP schemes previously recalled in Section \ref{sec:model}, respectively. 

\subsection{Solid body rotation}
\label{ssec:solid_body}

\begin{figure}[pos = H]
	\centering
	\includegraphics[width=0.9\textwidth]{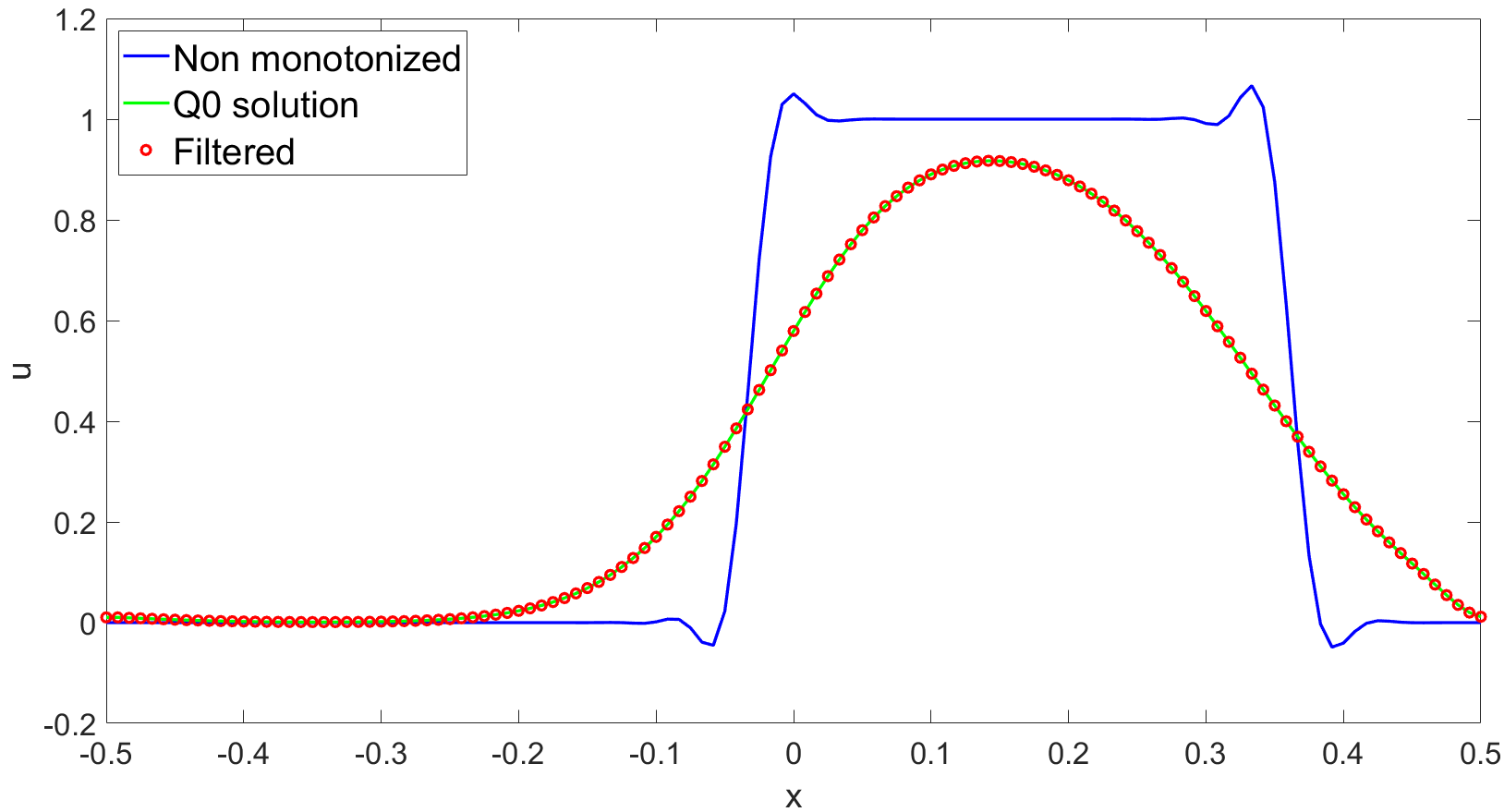} 
	\caption{Computational results for the solid body rotation at \(t = T_{f}\) with \(k = 1\) using filter \eqref{eq:filter} with \(\varepsilon = 5 h_{K}\). The red line denotes the non filtered \(Q_1\) solution, the green line denotes the \(Q_0\) solution, while the red dots represent the results of the simulation with the filtering approach.}
	\label{fig:Solid_Body_Rotation_epsilon_5}
\end{figure}

We consider a classical benchmark for convection schemes, the so-called solid body rotation, which has been studied in different configurations (see e.g. \cite{leveque:1996}, \cite{zalesak:1979}). A stationary velocity field is considered, representing a rotating flow with frequency \(\omega = \SI{1}{\per\second}\) around the point \(\left(0,0\right)\) on the domain \(\Omega = \left(-0.5,0.5\right)^2\). The initial datum is given by the following discontinuous function:
\[u_{0}(\mathbf{x}) = \begin{cases}
1 \qquad &\text{if } X^2 + Y^2 \le 1 \\
0 \qquad &\text{otherwise}
\end{cases}\]
where \(X = \frac{x - x_0}{\sigma}\) and \(Y = \frac{y - y_0}{\sigma}\) with \(x_0 = y_0 = \frac{1}{6}\) and \(\sigma = 0.2\). For this first test, the computational grid is composed by 120 elements along each direction with a time step such that the maximum Courant number is \(C \approx 0.1\). All the results are presented at \(T_{f} = \SI[parse-numbers=false]{2\pi}{\second}\), when one rotation has been completed, so that the solution coincides with the initial datum. We first consider the strategy \eqref{eq:filter} depicted in Section \ref{sec:model} with the filter function \(F_{1}(x)\) \eqref{eq:filter_function_a}, taking \(\varepsilon = 5 h_{K}\), as suggested in \cite{bokanowski:2016}, where \(h_{K} = \text{diam}(K)\). Figure \ref{fig:Solid_Body_Rotation_epsilon_5} compares the results at \(t = T_{f}\) of the filtering approach with the \(Q_1\) non monotonized solution and with the \(Q_0\) one. Recall that the finite element spaces are the ones defined in \eqref{eq:fe_spaces}. As one can easily notice, with this choice of the parameter, too much stabilization is added and therefore the filtered solution essentially coincides with the low order one.

\begin{figure}[pos = H]
	\centering
	\includegraphics[width=0.9\textwidth]{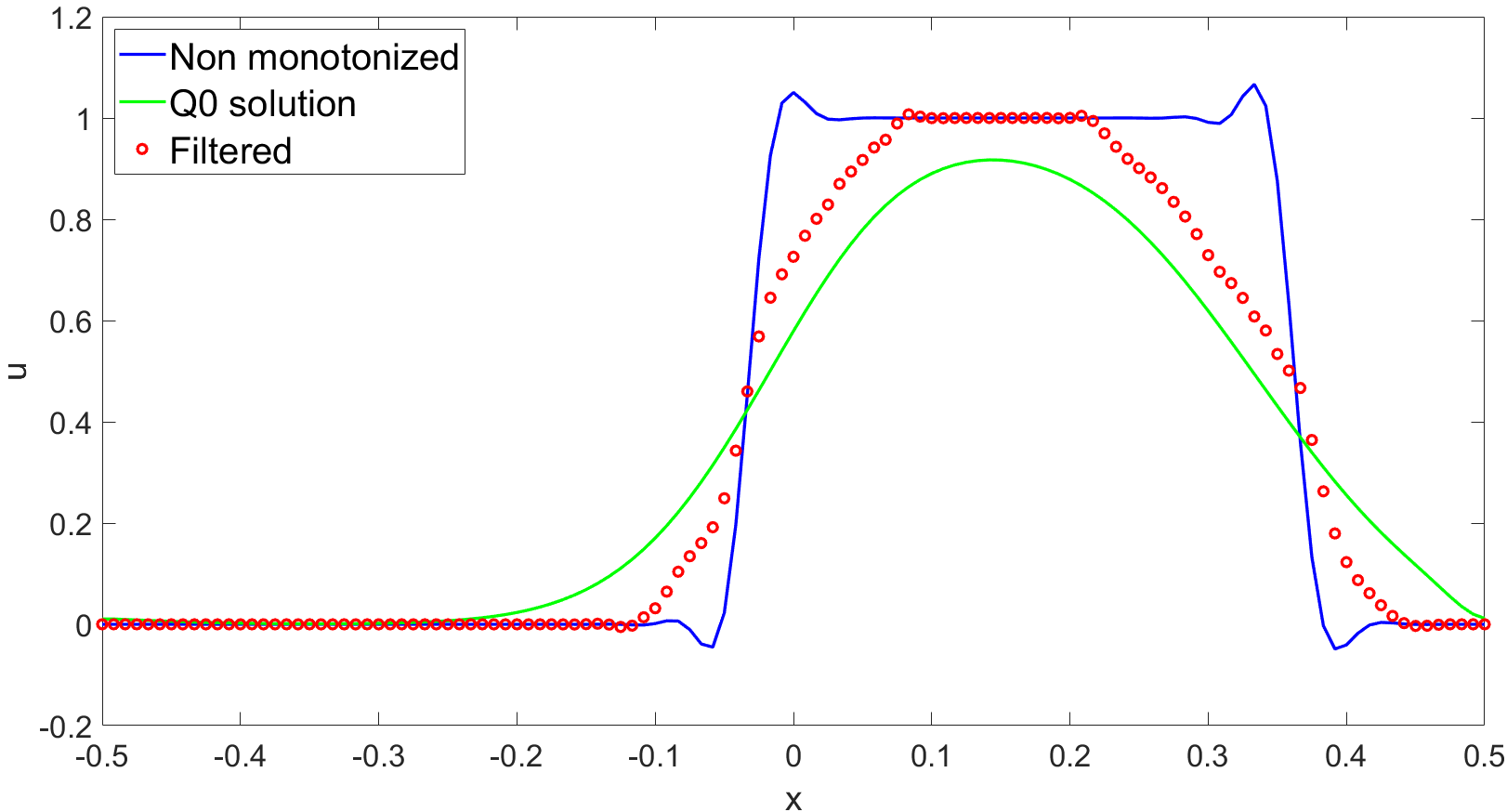} 
	\caption{Computational results for the solid body rotation at \(t = T_{f}\) with \(k = 1\) using filter \eqref{eq:filter} with \(\varepsilon = 10^{4} h_{K}\). The red line denotes the non filtered \(Q_1\) solution, the green line denotes the \(Q_0\) solution, while the red dots represent the results of the simulation with the filtering approach.}
	\label{fig:Solid_Body_Rotation_epsilon_10000}
\end{figure}

Increasing the value \(\varepsilon\) does not affect significantly the results until we take \(\varepsilon = 10^{4} h_{K}\): in this case, as evident from Figure \ref{fig:Solid_Body_Rotation_epsilon_10000}, the filtering approach works quite well since it is able to provide an essentially monotonic solution, as confirmed by Table \ref{tab:Solid_Body_Rotation_range_epsilon}, without smoothing it too much. 

\begin{figure}[pos = H]
	\centering
	\includegraphics[width=0.9\textwidth]{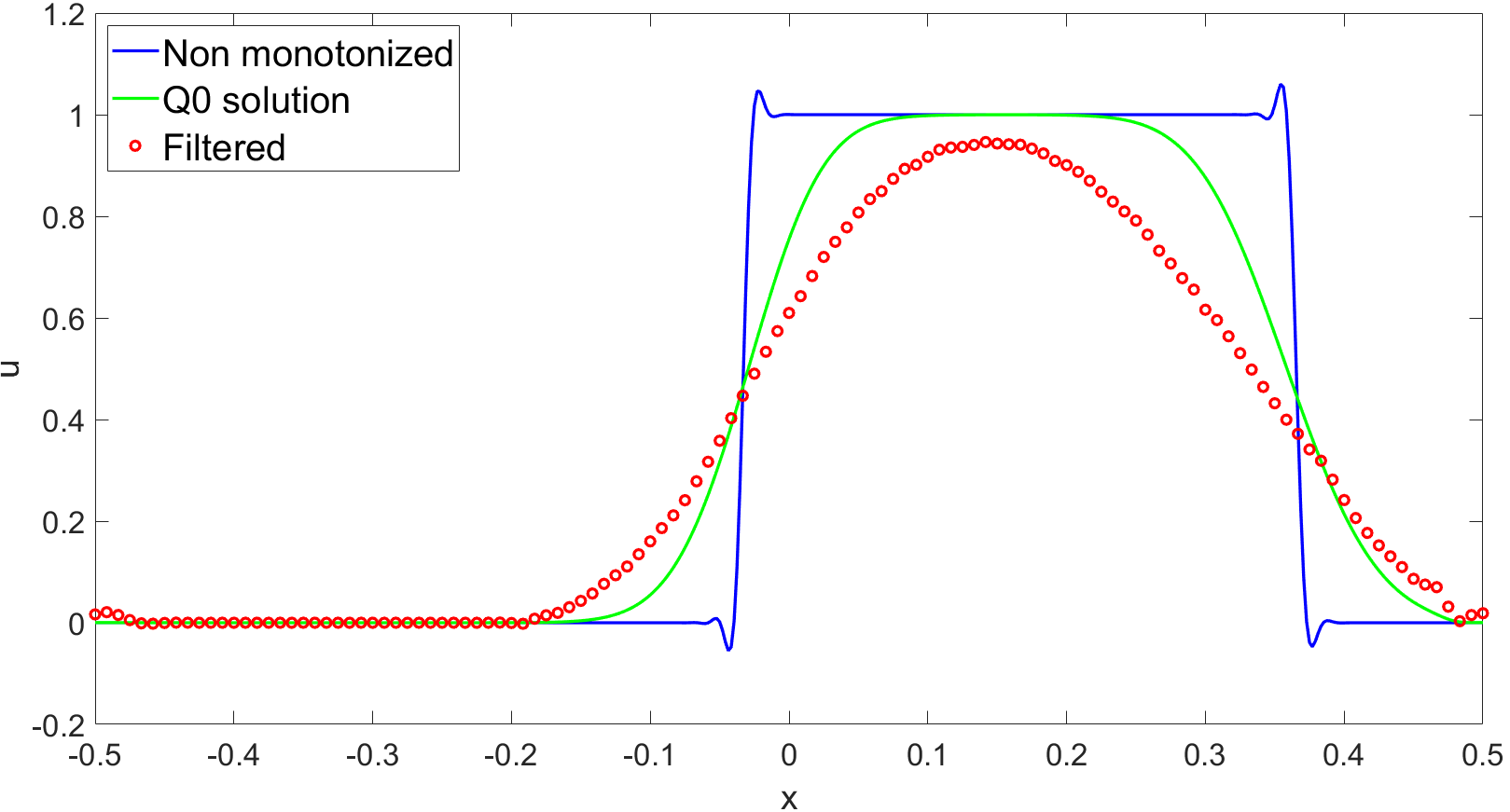} 
	\caption{Computational results for adaptive simulation of the solid body rotation at \(t = T_{f}\) with \(k = 1\) using filter \eqref{eq:filter} with \(\varepsilon = 10^{4} h_{K}\). The red line denotes the full resolution non filtered \(Q_1\) solution, the green line denotes the full resolution \(Q_0\) solution, while the red dots represent the results of the simulation with the filtering approach.}
	\label{fig:Solid_Body_Rotation_epsilon_10000_adaptive}
\end{figure}

The situation can be further improved using mesh adaptivity so as to start with a coarse mesh and perform refinement only in the zones where discontinuity is detected. The indicator is based on the gradient of the variable \(u\); more specifically, we define for each element \(K\)
\begin{equation}\label{eq:adaptive_criterion}
\eta_{K} = \max_{i \in \mathcal{N}_{K}} \left|\nabla u\right|_{i}
\end{equation}
as local refinement indicator, where \(\mathcal{N}_{K}\) denotes the set of nodes over the element \(K\). The initial mesh is composed by 120 elements along each direction and we allowed up to two local refinements. Figure \ref{fig:Solid_Body_Rotation_epsilon_10000_adaptive} shows that the results at \(t = T_{f}\) with a time step such that the maximum Courant number is \(C \approx 0.1\), using the value \(\varepsilon = 10^{4} h_{K}\) previously tested in the fixed grid configuration, compared with the full resolution \(Q_{1}\) non monotonized solution and the corresponding \(Q_{0}\) one. One can easily notice that in this specific configuration the value of \(\varepsilon\) is still too small and too much dissipation is provided.

\begin{figure}[pos = H]
	\centering
	\includegraphics[width=0.9\textwidth]{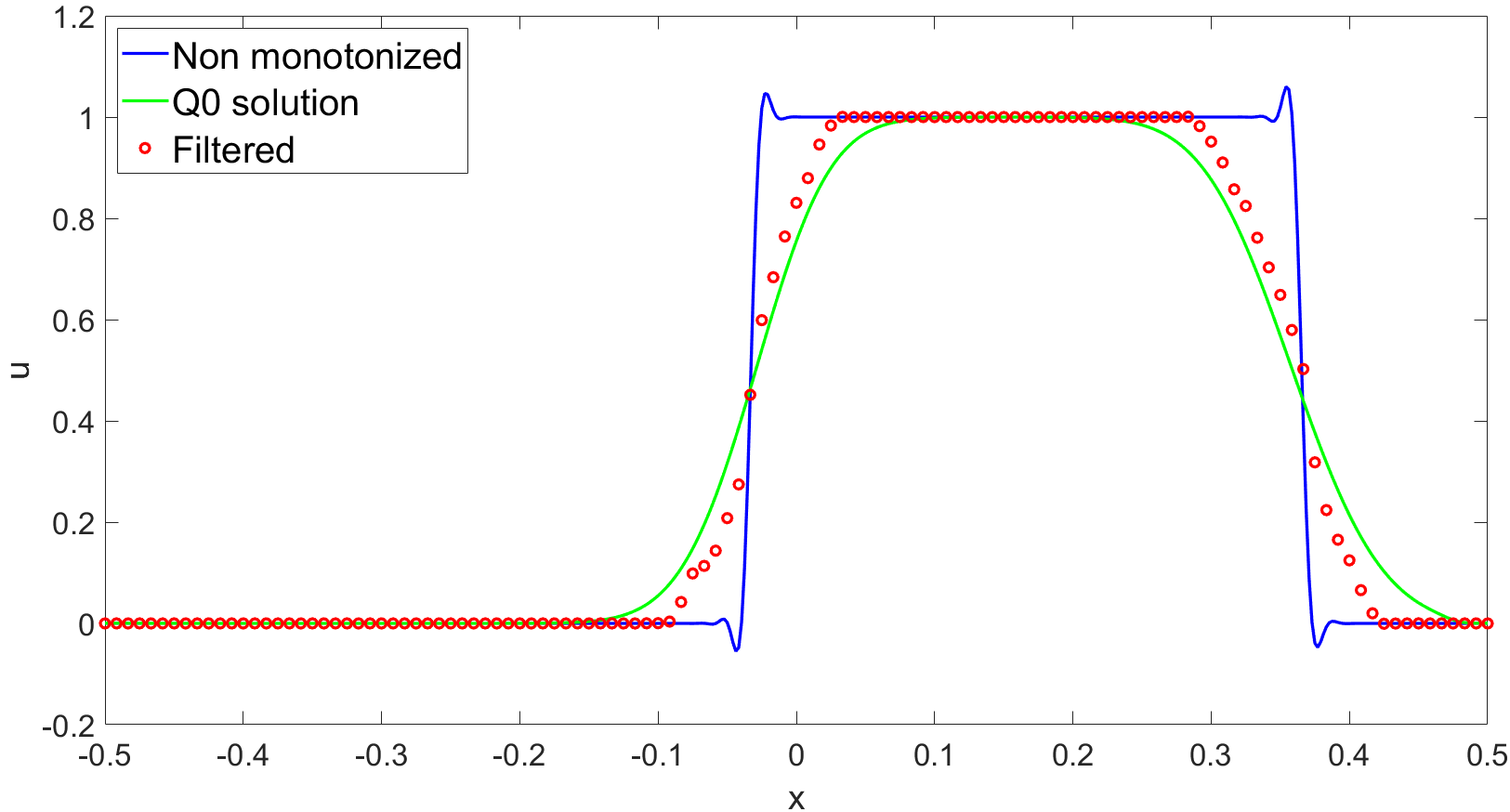} 
	\caption{Computational results for adaptive simulation of the solid body rotation at \(t = T_{f}\) using filter \eqref{eq:filter} with \(\varepsilon = 10^{5} h_{K}\). The red line denotes the full resolution non filtered \(Q_1\) solution, the green line denotes the full resolution \(Q_0\) solution, while the red dots represent the results of the simulation with the filtering approach.}
	\label{fig:Solid_Body_Rotation_epsilon_100000_adaptive}
\end{figure}

The situation improves increasing the value of \(\varepsilon\). Figure \ref{fig:Solid_Body_Rotation_epsilon_100000_adaptive} shows the results using \(\varepsilon = 10^{5} h_{K}\), where an essentially monotonic solution is retrieved. The values reported in Table \ref{tab:Solid_Body_Rotation_range_epsilon} confirm the better quality of the solution.

\begin{table}[pos = H]
	\centering
	\begin{tabular}{|c|c|c|}
		\hline
		\textbf{Value of } \(\varepsilon\) & \textbf{Maximum value of \(u\)}  & \textbf{Mininum value of \(u\)} \\
		\hline
		\(10^{4}h_{K}\) & \(1.0 + 7.3 \cdot 10^{-3}\) & \(0.0 - 5.2 \cdot 10^{-3}\)  \\
		\hline
		\(10^{5}h_{K}\) (adaptive) & \(1.0 + 1.0 \cdot 10^{-5}\) & \(0.0 - 1.9 \cdot 10^{-3}\)  \\
		\hline
	\end{tabular}
	\caption{Solid body rotation, maximum and minimum values for filtering approach \eqref{eq:filter} at \(t = T_{f}\) with \(k = 1\) both in case of fixed grid and adaptive simulations.}
	\label{tab:Solid_Body_Rotation_range_epsilon}
\end{table}

The very large value of \(\varepsilon\) that is necessary to achieve monotonicity suggests that the previous approach has shortcomings. We consider therefore the second strategy \eqref{eq:filter_rel} outlined in Section \ref{sec:model}. We start again from a fixed grid configuration, using the same mesh and the same time step previously described. After some sensitivity study, \(\beta = 0.4\) seems to yield an acceptable behaviour for the solution, as evident from Figure \ref{fig:Solid_Body_Rotation_beta0,4}. The discontinuity is less smeared out with respect to the \(Q_0\) solution, while avoiding the spurious oscillations and retrieving an essentially monotonic solution, as reported in Table \ref{tab:Solid_Body_Rotation_range_beta}. 

\begin{figure}[pos = H]
	\centering
	\includegraphics[width=0.9\textwidth]{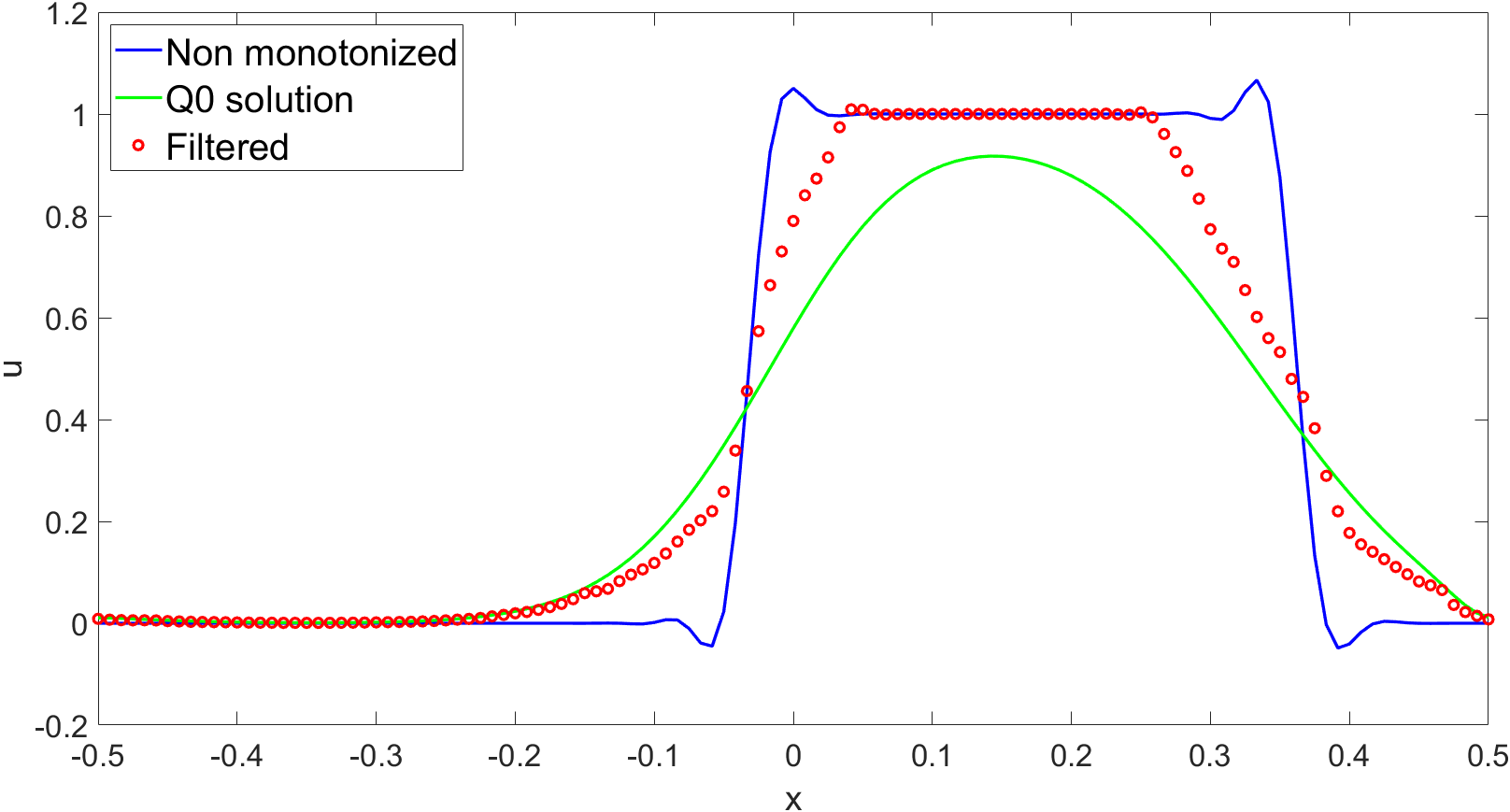} 
	\caption{Computational results for the solid body rotation at \(t = T_{f}\) with \(k = 1\) using filter \eqref{eq:filter_rel} with \(\beta = 0.4\). The red line denotes the non filtered \(Q_1\) solution, the green line denotes the \(Q_0\) solution, while the red dots represent the results of the simulation with the filtering approach.}
	\label{fig:Solid_Body_Rotation_beta0,4}
\end{figure}

Again, the \(h\)-adaptive version of the method, using the same configuration and the same refinement criterion previously described, provides better results, as confirmed by Table \ref{tab:Solid_Body_Rotation_range_beta}. The grid at \(t = T_{f}\) is reported in Figure \ref{fig:Solid_Body_Rotation_beta0,4_grid} and is composed by 28119 elements.

\begin{figure}[pos = H]
	\centering
	\includegraphics[width=0.9\textwidth]{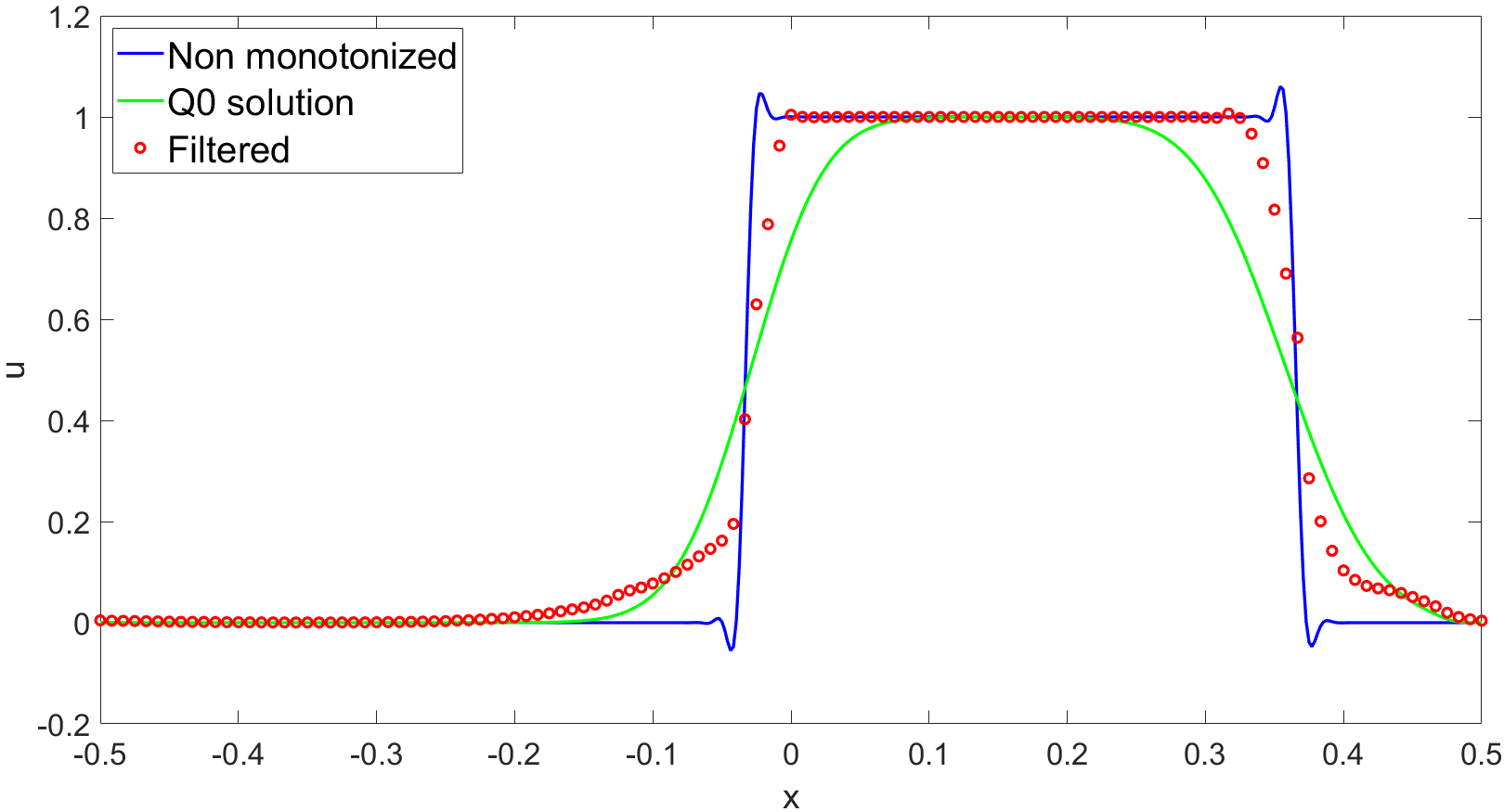} 
	\caption{Computational results for adaptive simulation of the solid body rotation at \(t = T_{f}\) with \(k = 1\) using filter \eqref{eq:filter_rel} with \(\beta = 0.4\). The red line denotes the full resolution non filtered \(Q_1\) solution, the green line denotes the full resolution \(Q_0\) solution, while the red dots represent the results of the simulation with the filtering approach.}
	\label{fig:Solid_Body_Rotation_beta0,4_adaptive}
\end{figure} 

\begin{table}[pos = H]
	\centering
	\begin{tabular}{|c|c|c|}
		\hline
		\textbf{Value of } \(\beta\) & \textbf{Maximum value of \(u\)}  & \textbf{Mininum value of \(u\)}   \\
		\hline
		\(0.4\) & \(1.0 + 9.2 \cdot 10^{-3}\) & \(0.0 + 7.6 \cdot 10^{-3}\)  \\
		\hline
		\(0.4\) (adaptive) & \(1.0 + 6.7 \cdot 10^{-3}\) & \(0.0 + 4.0 \cdot 10^{-4}\)  \\
		\hline
	\end{tabular}
	\caption{Solid body rotation, maximum and minimum values for filtering approach \eqref{eq:filter_rel} at \(t = T_{f}\) both in case of fixed grid and adaptive simulations.}
	\label{tab:Solid_Body_Rotation_range_beta}
\end{table}

\begin{figure}[pos = H]
	\centering
	\includegraphics[width=0.9\textwidth]{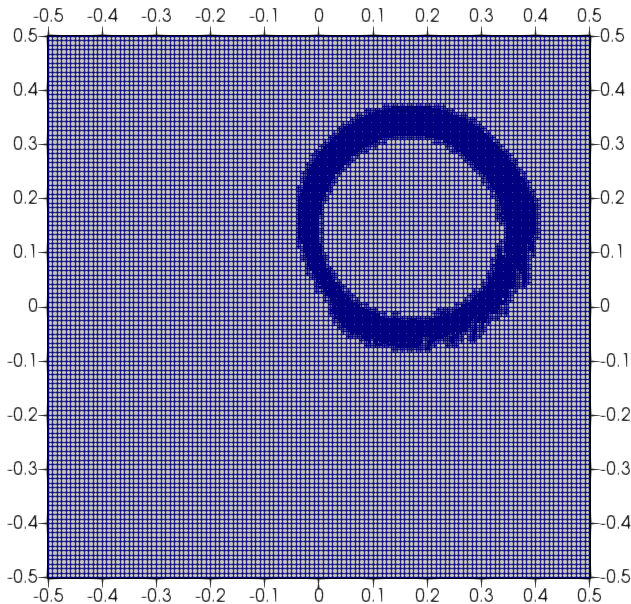} 
	\caption{Computational grid for adaptive simulation of the solid body rotation at \(t = T_{f}\) with \(k = 1\) using filter \eqref{eq:filter_rel} with \(\beta = 0.4\).}
	\label{fig:Solid_Body_Rotation_beta0,4_grid}
\end{figure} 

\begin{figure}[pos = H]
	\centering
	\includegraphics[width=0.9\textwidth]{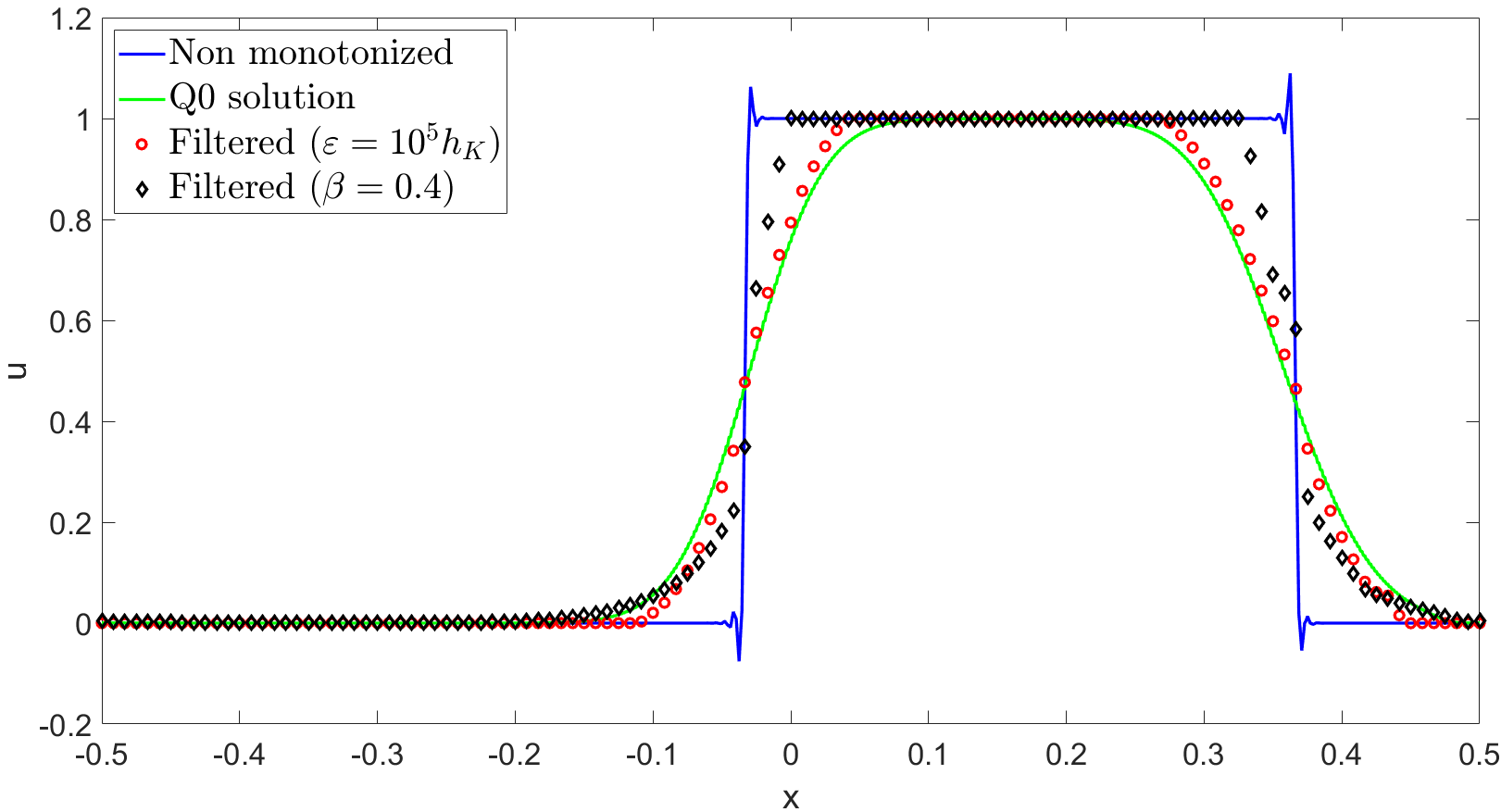} 
	\caption{Computational results for adaptive simulations of the solid body rotation at \(t = T_{f}\) with \(k = 2\). The red line denotes the full resolution non filtered \(Q_2\) solution, the green line denotes the full resolution \(Q_0\) solution, while the red dots and black diamonds represent the results of the simulation with the filtering approaches \eqref{eq:filter} and \eqref{eq:filter_rel} using \(\varepsilon = 10^{5} h_{K}\) and \(\beta = 0.4\), respectively.}
	\label{fig:Solid_Body_Rotation_comparison_Q2}
\end{figure}

The same test has been repeated using \(k = 2\), i.e. \(Q_{2}\) finite elements, and the third order SSP time discretization strategy briefly recalled in Section \ref{sec:model}. We present here only a comparison between the two strategies in case of an adaptive simulation using \(\varepsilon = 10^{5} h_{K}\) and \(\beta = 0.4\), respectively. Again, we started with a mesh composed by \(120\) elements along each direction, we allowed up to two local refinements and the employed time step is such that the maximum Courant number is \(C \approx 0.1\). Figure \ref{fig:Solid_Body_Rotation_comparison_Q2} shows the results with the two different approaches at \(t = T_{f}\), compared with a full resolution \(Q_{2}\) solution and the corresponding \(Q_{0}\) one. One can easily notice that both strategies provide an essentially monotonic result, as confirmed by Table \ref{tab:Solid_Body_Rotation_comparison_Q2}; moreover, the approach \eqref{eq:filter_rel} is characterized by a sharper transition zone and is therefore less dissipative, allowing to apply the filter on a reduced number of elements. Moreover, as evident from Table \ref{tab:Solid_Body_Rotation_comparison_Q2}, the filtering procedure \eqref{eq:filter_rel} appears to avoid undershoots and keeps a non negative solution, which is a crucial fact in many applications in order to preserve the physical meaning of the results. Hence, it will be the one used throughout the rest of the numerical experiments.

\begin{table}[pos = H]
	\centering
	\begin{tabular}{|c|c|c|}
		\hline
		\textbf{Value of the parameter} & \textbf{Maximum value of \(u\)}  & \textbf{Mininum value of \(u\)} \\
		\hline
		\(\beta = 0.4\) (adaptive)  &  \(1.0 + 1.5 \cdot 10^{-3}\) & \(0.0 + 3.8 \cdot 10^{-3}\)  \\
		\hline
		\(\varepsilon = 10^{5} h_{K}\) (adaptive)   &  \(1.0\) & \(0.0 - 1.9 \cdot 10^{-4}\)  \\
		\hline
	\end{tabular}
	\caption{Solid body rotation, maximum and minimum values for filtering approach \eqref{eq:filter} at \(t = T_{f}\) with \(k = 2\) in case of adaptive simulations.}
	\label{tab:Solid_Body_Rotation_comparison_Q2}
\end{table}

\begin{figure}[pos = H]
	\centering
	\includegraphics[width=0.9\textwidth]{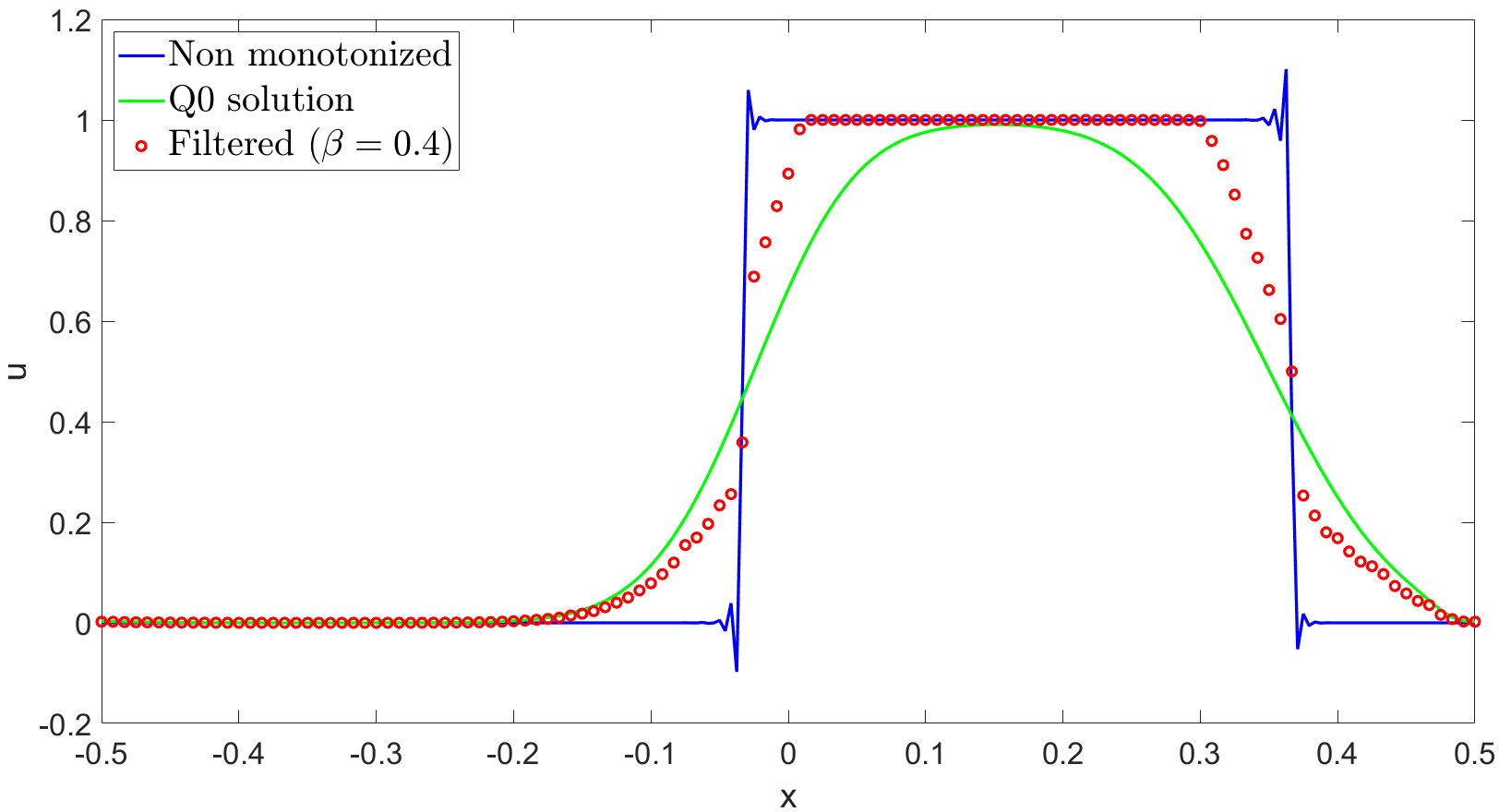} 
	\caption{Computational results for solid body rotation at \(t = T_{f}\) with \(k = 3\). The red line denotes the non filtered \(Q_3\) solution, the green line denotes the \(Q_0\) solution, while the red dots represent the results of the simulation with the filtering approaches \eqref{eq:filter_rel} using \(\beta = 0.4\).}
	\label{fig:Solid_Body_Rotation_Q3}
\end{figure}

For the sake of completeness, we report also the results obtained using the filtering approach \eqref{eq:filter_rel} with \(k = 3\), \(\beta = 0.4\) and a mesh composed by \(240\) elements along each direction. Figure \ref{fig:Solid_Body_Rotation_Q3} shows a comparison at \(t = T_{f}\) between the filtering solution, the non monotonized \(Q_{3}\) solution and the corresponding \(Q_{0}\) one. All the considerations made so far remain valid and the overshoots are further reduced with respect to Table \ref{tab:Solid_Body_Rotation_comparison_Q2}. The overhead with respect to the unfiltered DG scheme amounts to a factor \(\approx 1.5\) in terms of CPU time.

\begin{table}[pos = H]
	\centering
	\begin{tabular}{|c|c|c|}
		\hline
		\textbf{Value of the parameter} & \textbf{Maximum value of \(u\)}  & \textbf{Mininum value of \(u\)}   \\
		\hline
		\(\beta = 0.4\) & \(1.0 + 1.2 \cdot 10^{-4}\) & \(0.0 + 7.9 \cdot 10^{-6}\)  \\
		\hline
	\end{tabular}
	\caption{Solid body rotation, maximum and minimum values for filtering approach \eqref{eq:filter} at \(t = T_{f}\) with \(k = 3\) in case of adaptive simulations.}
	\label{tab:Solid_Body_Rotation_Q3}
\end{table}

\subsection{Smooth isentropic vortex}
\label{ssec:vortex}

The isentropic vortex problem is a classical benchmark for the two-dimensional compressible Euler equations introuced in \cite{shu:1998} (see also \cite{loubere:2014}) for which an analytic solution is available and can be therefore used to assess the convergence properties of a numerical scheme. The initial conditions are given as a perturbation of a reference state
\[\rho(\mathbf{x},0) = \rho_\infty + \delta\rho \qquad \mathbf{u}(\mathbf{x},0) = \mathbf{u}_{\infty} + \delta\mathbf{u} \qquad p(\mathbf{x},0) = p_{\infty} + \delta p.\]
The typical perturbation is defined as
\begin{equation}
\delta T = \frac{1 - \gamma}{8 \gamma \pi^2} \beta^2 e^{1 - r^2},
\end{equation} 
with \(r^2 = \left(x - x_0\right)^2 + \left(y - y_0\right)^2\) denoting the radial coordinate and \(\beta\) being the vortex strength. we set
\begin{equation}
\rho(\mathbf{x},0) = \left(1 + \delta T\right)^{\frac{1}{\gamma-1}} 
\ \ \ \
p(\mathbf{x},0) = \left(1 + \delta T\right)^{\frac{\gamma}{\gamma-1}}.
\end{equation}
For what concerns the velocity the typical perturbation is defined as
\begin{equation}
\delta \mathbf{u} = \beta \begin{pmatrix}
-\left(y- y_0\right) \\
\left(x - x_0\right)
\end{pmatrix}\frac{e^{\frac{1}{2}\left(1 - r^2\right)}}{2\pi},
\end{equation}
where \(x_0\) and \(y_0\) are the coordinates of the vortex centre. We consider the domain \(\Omega = \left(-5, 5\right)^2\) with periodic boundary conditions and we set \(\rho_\infty = 1\), \(p_\infty = 1\), \(x_0 = y_0 = 0\), \(\beta = 5\) and the final time \(T_f = 10\), so that the vortex is back to its original position. The simulations are performed at fixed Courant number \(C \approx 0.1\). Notice that, since we are in presence of a smooth solution, the problem should be simulated with effective high-order of accuracy. As evident from Tables \ref{tab:vortex_k1} and \ref{tab:vortex_k2} for the density, choosing a sufficient high value of the \(\beta\) parameters avoids to activate the filter and allows hence to achieve the expected convergence rates, whereas, for too small values, the overall convergence rates are affected by the \(Q_0\) solution. Analogous results are obtained for the momentum and for the energy. The results compare well with the one reported in \cite{dumbser:2014} and \cite{zanotti:2015}.

\begin{table}[pos = H]
	\centering
	\begin{tabular}{|c|c|c|c|c|c|c|c|c|c|}
		\hline
		\(\beta_\rho\) & \(\beta_{\rho\mathbf{u}}\) & \(\beta_{\rho E}\) & \(N_{el}\) & \(L^1\) rel. error \(\rho\) & \(L^1\) rate \(\rho\) & \(L^2\) rel. error \(\rho\) & \(L^2\) rate \(\rho\) & \(L^\infty\) rel. error \(\rho\) & \(L^\infty\) rate \(\rho\) \\
		\hline
		& & & \(20\) & \(8.44 \cdot 10^{-3}\) &  & \(1.99 \cdot 10^{-2}\) & & \(1.53 \cdot 10^{-1}\) & \\
		0.7 & 0.7 & 0.7 & \(40\) & \(6.23 \cdot 10^{-3}\) & \(0.44\) & \({1.57} \cdot 10^{-2}\) & \(0.34\) & \(1.33 \cdot 10^{-1}\) & \(0.20\) \\
		& & & \(80\) & \(5.04 \cdot 10^{-3}\) & \(0.31\) & \(1.25 \cdot 10^{-2}\) & \(0.33\) & \(1.09 \cdot 10^{-1}\) & \(0.29\) \\
		& & & \(160\) & \(2.73 \cdot 10^{-3}\) & \(0.88\) & \(6.60 \cdot 10^{-3}\) & \(0.92\) & \(5.92 \cdot 10^{-2}\) & \(0.88\) \\
		\hline
		& & & \(20\) & \(3.63 \cdot 10^{-3}\) &  & \(8.62 \cdot 10^{-3}\) & & \(6.31 \cdot 10^{-2}\) & \\
		1.0 & 1.0 & 1.0 & \(40\) & \(8.02 \cdot 10^{-4}\) & \(2.18\) & \({1.81} \cdot 10^{-3}\) & \(2.25\) & \(1.26 \cdot 10^{-2}\) & \(2.32\) \\
		& & & \(80\) & \(1.79 \cdot 10^{-4}\) & \(2.16\) & \(3.89 \cdot 10^{-4}\) & \(2.22\) & \(2.76 \cdot 10^{-3}\) & \(2.19\) \\
		& & & \(160\) & \(4.19 \cdot 10^{-5}\) & \(2.09\) & \(8.94 \cdot 10^{-5}\) & \(2.12\) & \(7.00 \cdot 10^{-4}\) & \(1.98\) \\
		\hline  
	\end{tabular}
	\caption{Convergence test for the isentropic vortex at \(C \approx 0.1\) with  \(k = 1\). Relative errors for the density in \(L^{1}, L^{2}\) and \(L^{\infty}\) norm. \(N_{el}\) denotes the number of elements along each direction.}
	\label{tab:vortex_k1}
\end{table}
\begin{table}[pos = H]
	\centering
	\begin{tabular}{|c|c|c|c|c|c|c|c|c|c|}
		\hline
		\(\beta_\rho\) & \(\beta_{\rho\mathbf{u}}\) & \(\beta_{\rho E}\) & \(N_{el}\) & \(L^1\) rel. error \(\rho\) & \(L^1\) rate \(\rho\) & \(L^2\) rel. error \(\rho\) & \(L^2\) rate \(\rho\) & \(L^\infty\) rel. error \(\rho\) & \(L^\infty\) rate \(\rho\) \\
		\hline
		& & & \(20\) & \(7.01 \cdot 10^{-3}\) &  & \(1.77 \cdot 10^{-2}\) & & \(1.63 \cdot 10^{-1}\) & \\
		0.7 & 0.7 & 0.7 & \(40\) & \(6.32 \cdot 10^{-3}\) & \(0.15\) & \({1.61} \cdot 10^{-2}\) & \(0.14\) & \(1.45 \cdot 10^{-1}\) & \(0.17\) \\
		& & & \(80\) & \(5.10 \cdot 10^{-3}\) & \(0.31\) & \(1.27 \cdot 10^{-2}\) & \(0.34\) & \(1.14 \cdot 10^{-1}\) & \(0.35\) \\
		& & & \(160\) & \(2.75 \cdot 10^{-3}\) & \(0.89\) & \(6.64 \cdot 10^{-3}\) & \(0.94\) & \(6.07 \cdot 10^{-2}\) & \(0.91\) \\
		\hline
		& & & \(20\) & \(2.12 \cdot 10^{-4}\) &  & \(3.86 \cdot 10^{-4}\) & & \(3.87 \cdot 10^{-3}\) & \\
		1.0 & 1.0 & 1.0 & \(40\) & \(3.03 \cdot 10^{-5}\) & \(2.82\) & \({5.89} \cdot 10^{-5}\) & \(2.71\) & \(6.81 \cdot 10^{-4}\) & \(2.51\) \\
		& & & \(80\) & \(4.57 \cdot 10^{-6}\) & \(2.73\) & \(1.07 \cdot 10^{-5}\) & \(2.45\) & \(1.27 \cdot 10^{-4}\) & \(2.42\) \\
		& & & \(160\) & \(9.22 \cdot 10^{-7}\) & \(2.31\) & \(1.85 \cdot 10^{-6}\) & \(2.53\) & \(2.19 \cdot 10^{-5}\) & \(2.54\) \\
		& & & \(320\) & \(1.45 \cdot 10^{-7}\) & \(2.67\) & \(2.88 \cdot 10^{-7}\) & \(2.68\) & \(4.08 \cdot 10^{-6}\) & \(2.43\) \\
		\hline  
	\end{tabular}
	\caption{Convergence test for the isentropic vortex at \(C \approx 0.1\) with \(k = 2\). Relative errors for the density in \(L^{1}, L^{2}\) and \(L^{\infty}\) norm. \(N_{el}\) denotes the number of elements along each direction.}
	\label{tab:vortex_k2}
\end{table}

\subsection{Sod shock tube problem}
\label{ssec:sod}

\begin{figure}[pos = H]
	\centering
	\includegraphics[width=0.9\textwidth]{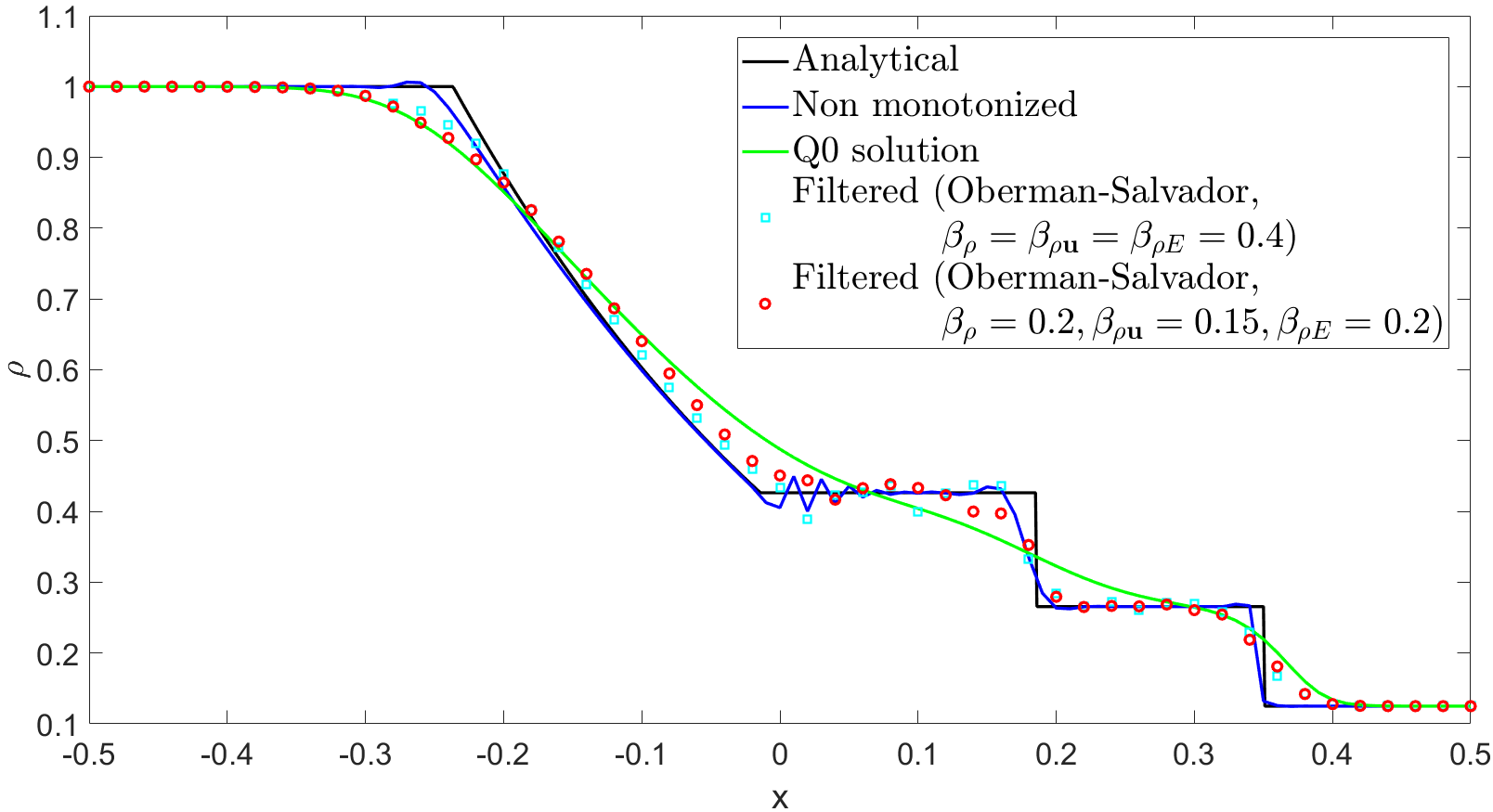} 
	\caption{Computational results for Sod shock tube problem at \(t = T_{f}\) with \(k = 1\). The black line reports the analytical solution, the red line denotes the non filtered \(Q_1\) solution, the green line denotes the \(Q_0\) solution, the cyan squares represent the results of the simulation with the filtering approach \eqref{eq:filter_function_a} using \(\beta_{\rho} = \beta_{\rho\mathbf{u}} = \beta_{\rho E} = 0.4\), whereas the red circles represent the results of the simulation with the filtering approach \eqref{eq:filter_function_a} using \(\beta_{\rho} = 0.2, \beta_{\rho\mathbf{u}} = 0.15, \beta_{\rho E} = 0.2\).}
	\label{fig:Sod_ObermanSalvador}
\end{figure}

We consider now the classical Sod shock tube problem proposed by \cite{sod:1978} in order to assess the capability of the filtering approach to reproduce correctly 1D waves such as shocks, contact discontinuities or rarefaction waves. It consists of a right-moving shock wave, an intermediate contact discontinuity and a left-moving rarefaction fan. The computational domain is \(\Omega = \left(-0.5, 0.5\right)\), the final time is \(T_{f} = \SI{0.2}{\second}\) and the initial condition is given as follows:
\begin{equation}
\left(\rho_{0}, u_{0}, p_{0}\right) = \begin{cases}
\left(1, 0, 1\right) \qquad &\text{if } x < 0 \\
\left(0.125, 0, 1\right) \qquad &\text{if } x > 0,
\end{cases}
\end{equation} 
in terms of density, velocity and pressure, respectively. Dirichlet boundary conditions are imposed. We use as numerical flux the Rusanov \cite{rusanov:1962} flux. We start with a mesh composed by \(100\) elements and a time-step equal to \(5 \cdot 10^{-4} \hspace{0.1cm} \SI{}{\second}\) and \(k = 1\), yielding a maximum Courant number \(C \approx 0.12\). Figure \ref{fig:Sod_ObermanSalvador} shows the results at \(t = T_{f}\) for the density of a simulation using \(\beta_{\rho} = \beta_{\rho\mathbf{u}} = \beta_{\rho E} = 0.4\). One can easily notice the presence of significant under- and over-shoots. This suggests that we need to decrease the value of the parameter \(\beta_{\rho}\) in order to achieve a monotonic solution. The same considerations hold also for the velocity and the pressure. After some sensitivity study, the combination \(\beta_{\rho} = 0.2, \beta_{\rho\mathbf{u}} = 0.15, \beta_{\rho E} = 0.2\) could be shown to provide a better quality solution with significantly reduced under- and over-shoots, as reported in Figure \ref{fig:Sod_ObermanSalvador}.

\FloatBarrier

\begin{figure}[pos = H]
	\begin{subfigure}[]{\textwidth}
		\centering
		\includegraphics[width=0.75\textwidth]{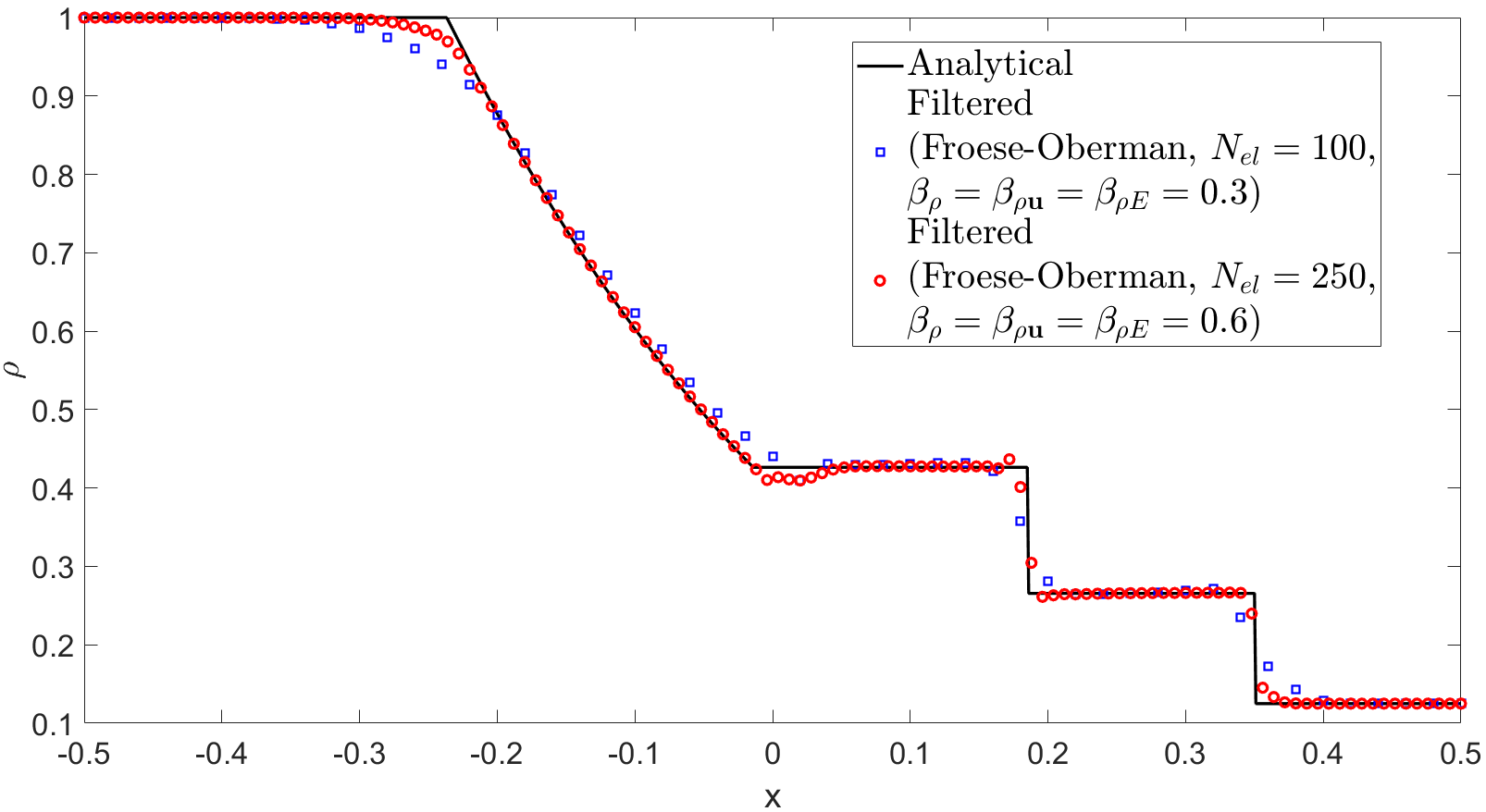} a)
	\end{subfigure}
	\begin{subfigure}[]{\textwidth}
		\centering
		\includegraphics[width=0.75\textwidth]{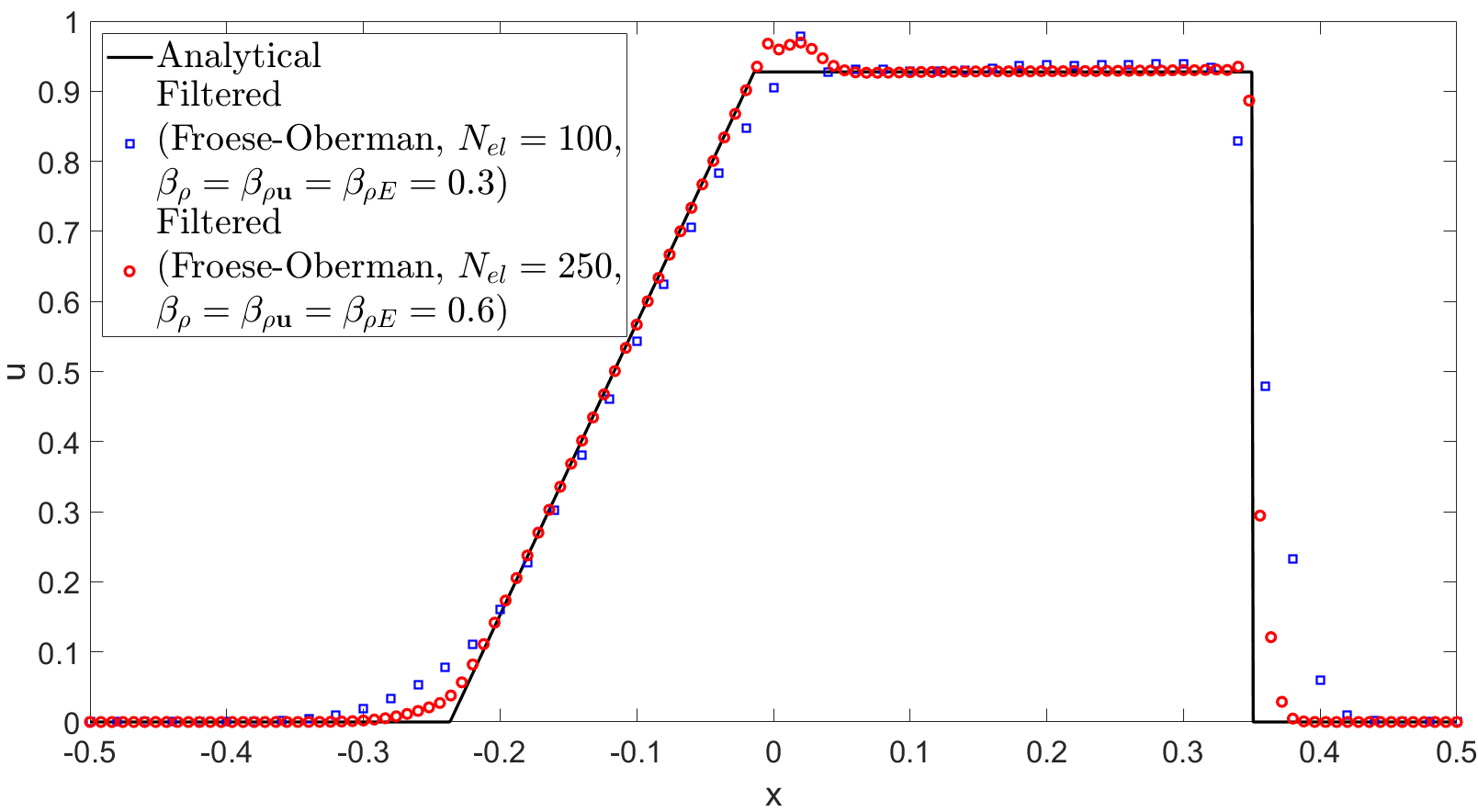} b)
	\end{subfigure} 
	\begin{subfigure}[]{\textwidth}
		\centering
		\includegraphics[width=0.75\textwidth]{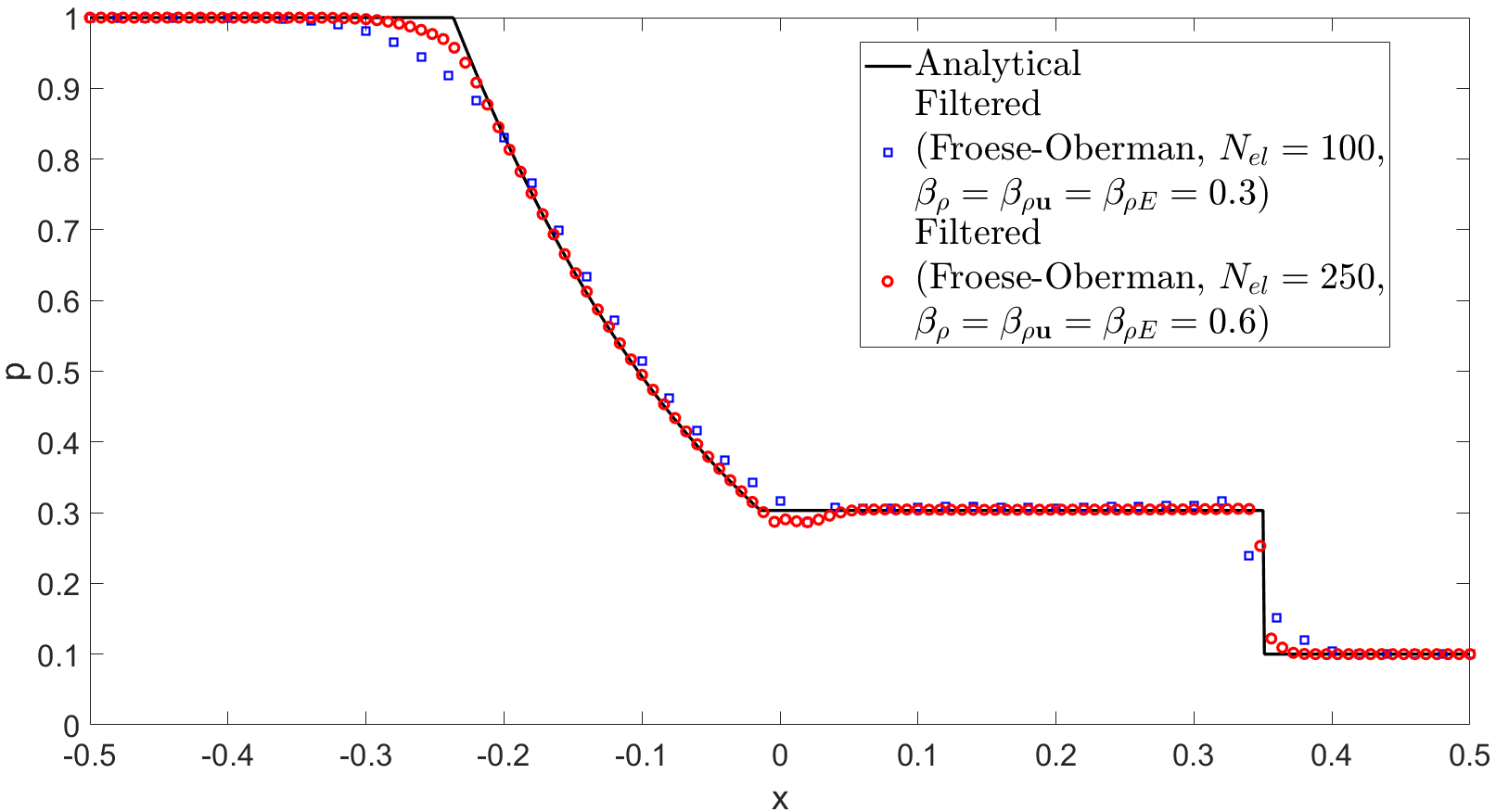} c)
	\end{subfigure}
	\caption{Computational results for Sod shock tube problem at \(t = T_{f}\) with \(k = 1\), a) density, b) velocity, c) pressure. The black line reports the analytical solution, the red line denotes the non filtered \(Q_1\) solution, the blue dots are the results employing Froese and Oberman's filter function \eqref{eq:filter_function_b} and using \(\beta_{\rho} = 0.3, \beta_{\rho\mathbf{u}} = 0.3\) and \(\beta_{\rho E} = 0.3\) with 100 elements, while the red dots represent the results employing Froese and Oberman's filter function \eqref{eq:filter_function_b} and using \(\beta_{\rho} = 0.6, \beta_{\rho\mathbf{u}} = 0.6\) and \(\beta_{\rho E} = 0.6\) with 250 elements.}
	\label{fig:Sod_FO}
\end{figure}

The situation can be further improved employing the Froese and Oberman's filter function \(F_2(x)\), that is continuous and provides therefore a smoother transition between the high order and the low order solutions. This allows also to increase the values of the parameters \(\beta_{\rho}, \beta_{\rho\mathbf{u}} \) and \(\beta_{\rho E}\). Figure \ref{fig:Sod_FO} shows the results at \(t = T_{f}\) using \(\beta_{\rho} = \beta_{\rho\mathbf{u}} = \beta_{\rho E} = 0.3\) and one can easily notice that the shock wave and the contact discontinuity are resolved in a sharper manner with only slight undershoots for density and pressure and overshoots for the velocity in the tail of the rarefaction fan. Table \ref{tab:Sod_beta0,3_beta0,3_beta0,3} reports the maximum and the minimum values for density, velocity and pressure, as well as the \(L^{\infty}\) norm errors, which confirm the good results of the proposed method, also in comparison with the results obtained in \cite{loubere:2014} with the classical ADER-MOOD and ADER-WENO schemes. Figure \ref{fig:Sod_FO} reports also the results at \(t = T_{f}\) using \(250\) elements, a time step equal to \(2 \cdot 10^{-4} \hspace{0.1cm} \SI{}{\second}\) and the following parameters: \(\beta_{\rho} = 0.6, \beta_{\rho\mathbf{u}} = 0.6, \beta_{\rho E} = 0.6\). It can be easily noticed that, as expected by increasing the resolution, the discontinuities are better retrieved. The values reported in Table \ref{tab:Sod_beta0,6_beta0,6_beta0,6} confirm the improved results.
  
\begin{table}[pos = H]
	\centering
	\begin{tabular}{|c|c|c|c|c|c|}
		\hline
		\textbf{Variable} & \textbf{Maximum value} & \textbf{Mininum value} & \(\mathbf{L^{\infty}}\) \textbf{error} & \(\mathbf{L^{\infty}}\) \textbf{error ADER-MOOD} & \(\mathbf{L^{\infty}}\) \textbf{error ADER-WENO} \\
		\hline
		\(\rho\) & \(1.0\) & \(0.125\) & \(6.9 \cdot 10^{-2}\) & \(1.1 \cdot 10^{-1}\) & \(1.2 \cdot 10^{-1}\)\\
		\hline
		\(u\) & \(0.9275 + 5.0 \cdot 10^{-2}\) & \(0.0\) &  \(4.8 \cdot 10^{-1}\) &  & \\
		\hline
		\(p\) & \(1.0\) & \(0.1\) & \(8.2 \cdot 10^{-2}\) & & \\
		\hline
	\end{tabular}
	\caption{Computational results for Sod shock tube problem at \(t = T_{f}\) with \(k = 1\) employing Froese and Oberman's filter function and using \(\beta_{\rho} = 0.3\), \(\beta_{\rho\mathbf{u}} = 0.3\) and \(\beta_{\rho E} = 0.3\). Results for ADER-WENO and ADER-MOOD schemes from \cite{loubere:2014}. Maximum and minimum values are referred to the corresponding variable declared in the left column.}
	\label{tab:Sod_beta0,3_beta0,3_beta0,3}
\end{table}

\begin{table}[pos = H]
	\centering
	\begin{tabular}{|c|c|c|c|}
		\hline
		\textbf{Variable} & \textbf{Maximum value}  & \textbf{Mininum value} & \(\mathbf{L^{\infty}}\) \textbf{error}  \\
		\hline
		\(\rho\) & \(1.0\) & \(0.125\) & \(3.9 \cdot 10^{-2}\) \\
		\hline
		\(u\) & \(0.9275 + 4.2 \cdot 10^{-2}\) & \(0.0\) &  \(2.9 \cdot 10^{-1}\)   \\
		\hline
		\(p\) & \(1.0\) & \(0.1\) & \(5.0 \cdot 10^{-2}\) \\
		\hline
	\end{tabular}
	\caption{Computational results for Sod shock tube problem at \(t = T_{f}\) with \(k = 1\) and \(250\) elements employing Froese and Oberman's filter function and using \(\beta_{\rho} = 0.6\), \(\beta_{\rho\mathbf{u}} = 0.6\) and \(\beta_{\rho E} = 0.6\). Maximum and minimum values are referred to the corresponding variable declared in the left column.}
	\label{tab:Sod_beta0,6_beta0,6_beta0,6}
\end{table}

The same test has been repeated using \(k = 2\) and the third order SSP time discretization scheme. Figure \ref{fig:Sod_beta1,4_beta1,4_beta1,4} reports the results at \(t = T_{f}\) using \(250\) elements, a time step equal to \(1 \cdot 10^{-4} \hspace{0.1cm} \SI{}{\second}\) and \(\beta_{\rho} = \beta_{\rho\mathbf{u}} = \beta_{\rho E} = 1.4\). The under- and over-shoots are significantly reduced and a good agreement with the analytical solution is established. The larger values of \(\beta\) parameters can be explained by considering that the increase of the polynomial degree leads generally to a more accurate solution with relatively large under- and over-shoots localized in a narrow region, where the low order solution has to be considered. Approximately, the filtering is applied on the 10\% of the degrees of freedom and the overhead with respect to the non monotonized scheme corresponds to a factor \(\approx 1.25\) in terms of CPU time. Both data compare quite well with the one reported in \cite{zanotti:2015} for the ADER-WENO approach, where the 15 \% of the cells was limited.

\begin{table}[pos = H]
	\centering
	\begin{tabular}{|c|c|c|c|}
		\hline
		\textbf{Variable} & \textbf{Maximum value}  & \textbf{Mininum value} & \(\mathbf{L^{\infty}}\) \textbf{error}  \\
		\hline
		\(\rho\)    &  \(1.0 + 2.0 \cdot 10^{-3}\) & \(0.125 - 2.0 \cdot 10^{-4}\) & \(1.6 \cdot 10^{-2}\) \\
		\hline
		\(u\)   &  \(0.9275 + 8.0 \cdot 10^{-3}\) & \(0.0 - 2.3 \cdot 10^{-3}\) &  \(2.2 \cdot 10^{-2}\)   \\
		\hline
		\(p\)    &  \(1.0 + 2.8 \cdot 10^{-3}\) & \(0.1 - 2.0 \cdot 10^{-4}\) & \(1.3 \cdot 10^{-2}\) \\
		\hline
	\end{tabular}
	\caption{Computational results for Sod shock tube problem at \(t = T_{f}\) with 250 elements and \(k = 2\) employing Froese and Oberman's filter function and using \(\beta_{\rho} = 1.4\), \(\beta_{\rho\mathbf{u}} = 1.4\) and \(\beta_{\rho E} = 1.4\). Maximum and minimum values are referred to the corresponding variable declared in the left column.}
	\label{tab:Sod_beta1,4_beta1,4_beta1,4}
\end{table}
\FloatBarrier

\begin{figure}[pos = H]
	\begin{subfigure}[]{\textwidth}
		\centering
		\includegraphics[width=0.75\textwidth]{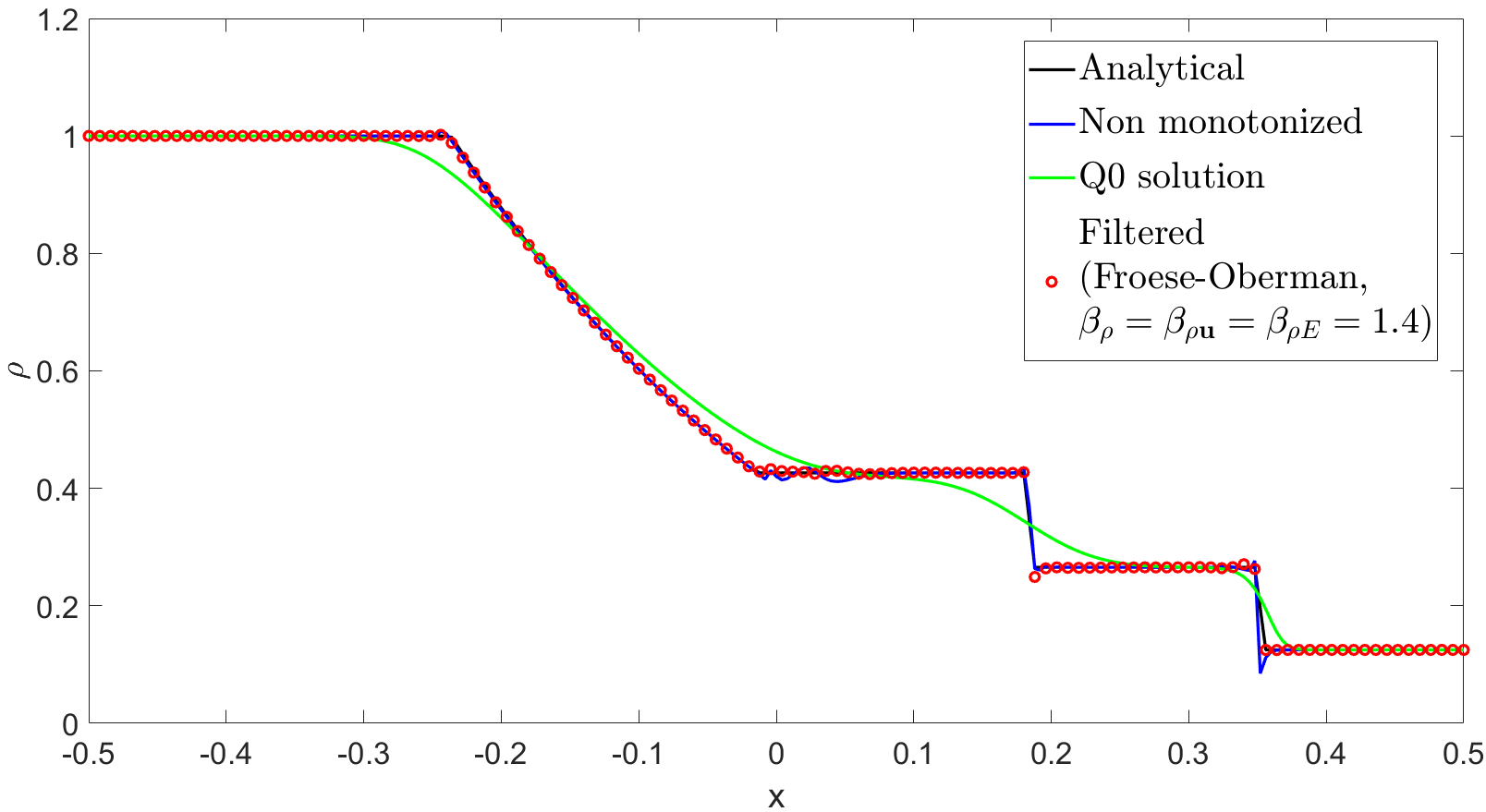} a) 
	\end{subfigure}
	\begin{subfigure}[]{\textwidth}
		\centering
		\includegraphics[width=0.75\textwidth]{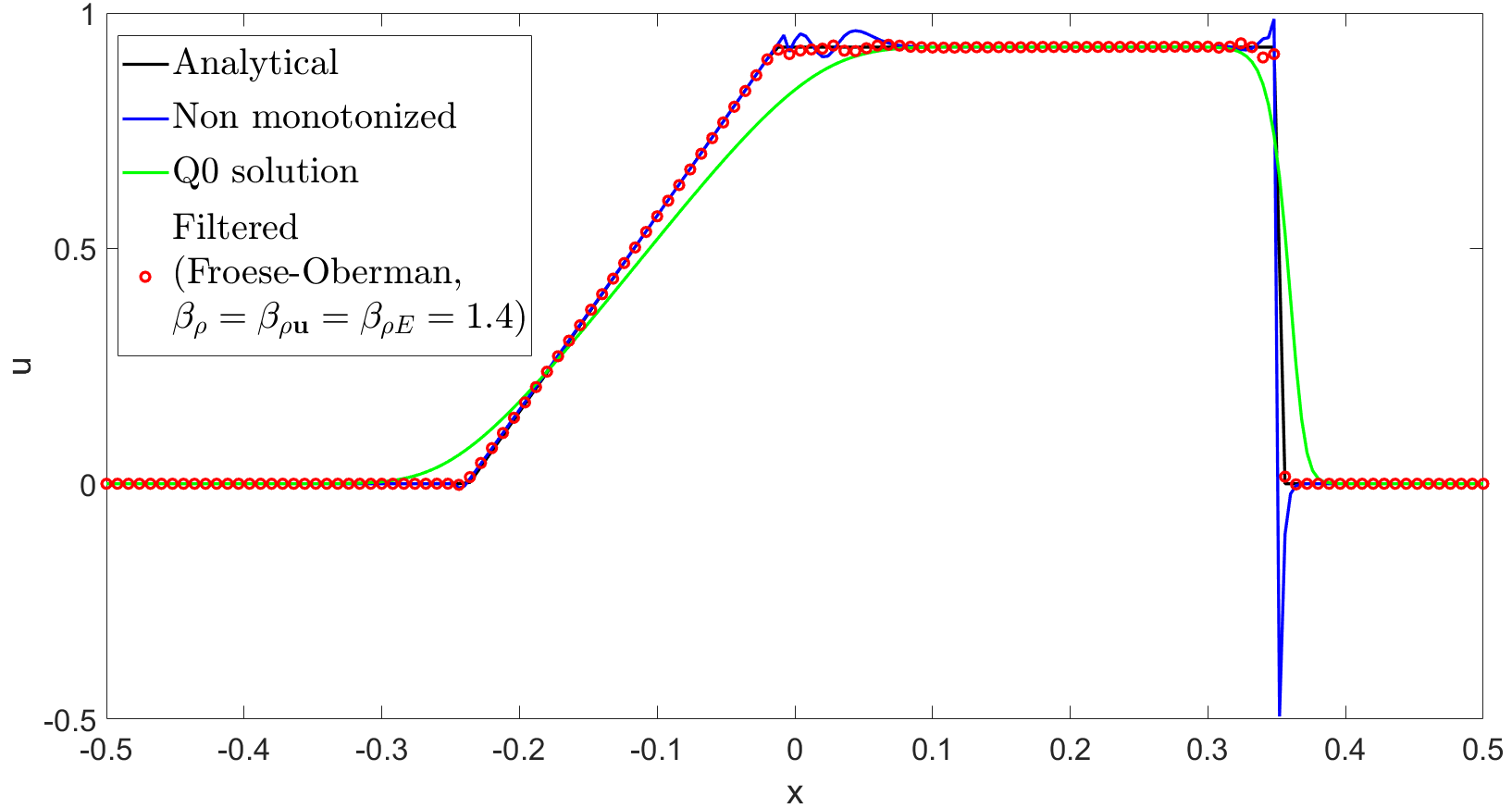} b)
	\end{subfigure}
	\begin{subfigure}[]{\textwidth}
		\centering
		\includegraphics[width=0.75\textwidth]{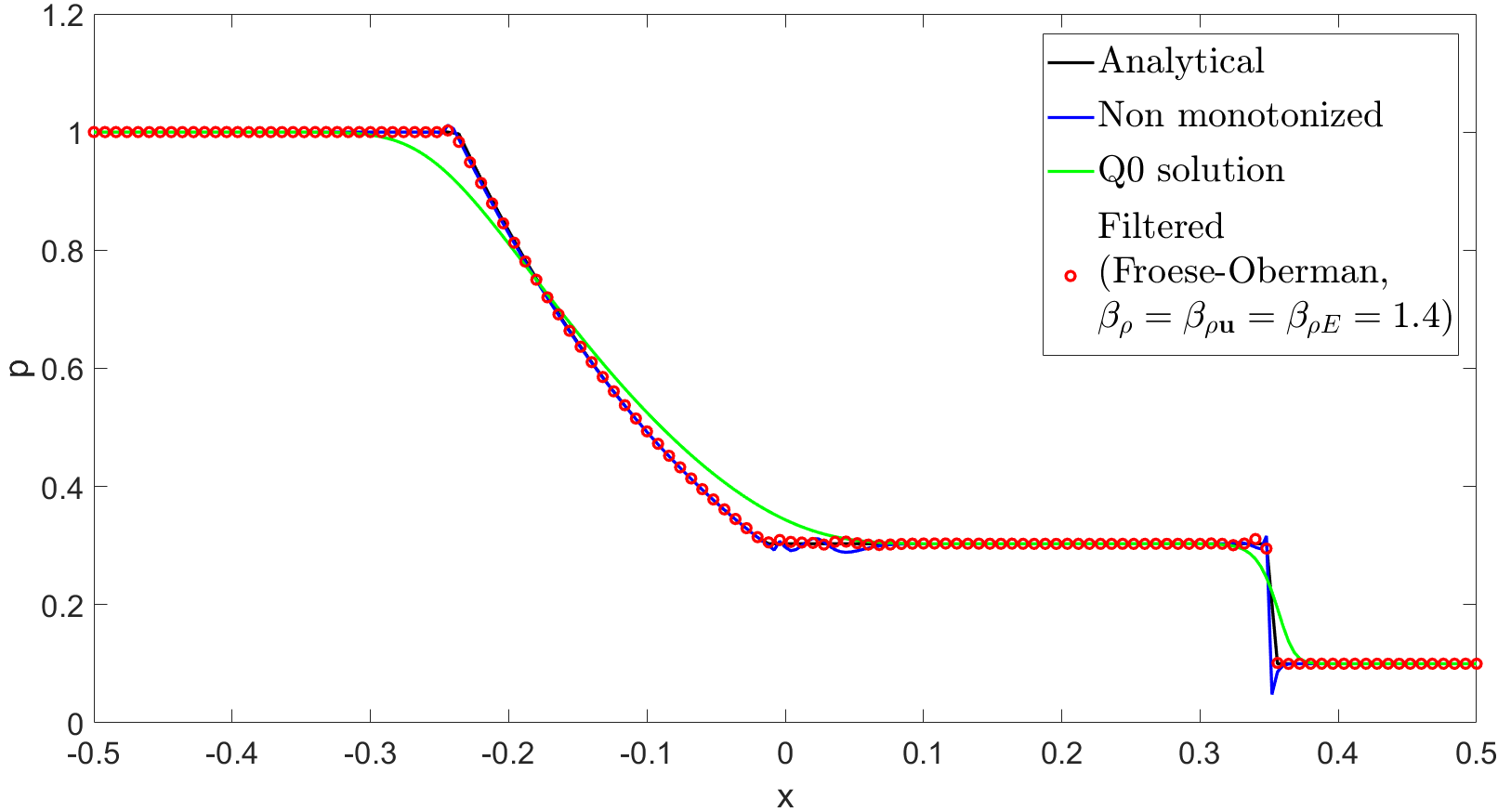} c)
	\end{subfigure}
	\caption{Computational results for Sod shock tube problem at \(t = T_{f}\) with \(k = 2\), a) density, b) velocity, c) pressure. The black line reports the analytical solution, the red line denotes the non filtered \(Q_2\) solution, the green line denotes the \(Q_0\) solution, while the red dots represent the results of the simulation with Froese and Oberman's filter function using \(\beta_{\rho} = 1.4\), \(\beta_{\rho\mathbf{u}} = 1.4\) and \(\beta_{\rho E} = 1.4\).}
	\label{fig:Sod_beta1,4_beta1,4_beta1,4}
\end{figure}

\subsection{Circular explosion problem}
\label{ssec:explosion}

\begin{figure}[pos = H]
	\begin{subfigure}[]{\textwidth}
		\centering
		\includegraphics[width=0.7\textwidth]{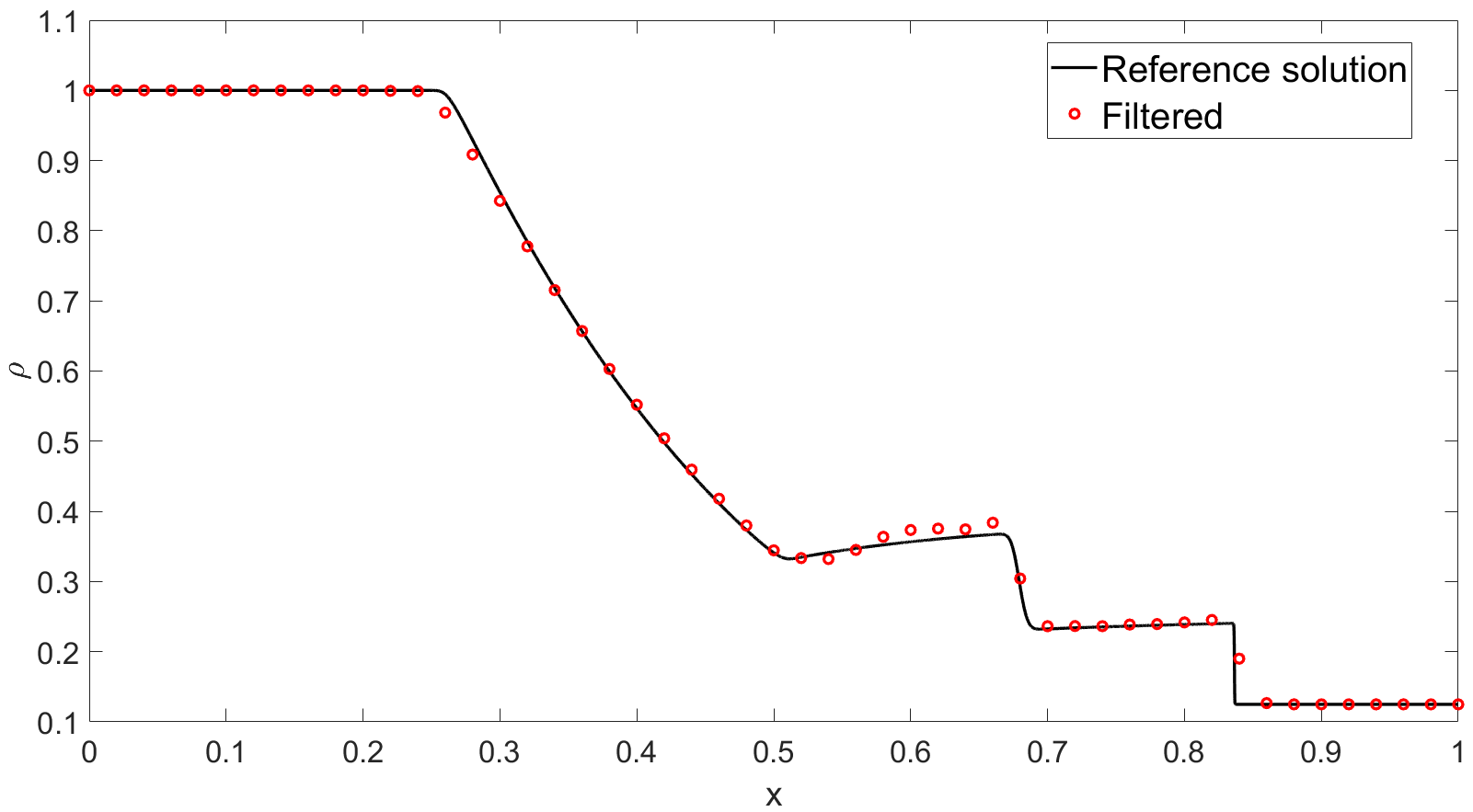} a)
	\end{subfigure}	
	\begin{subfigure}[]{\textwidth}
		\centering 
		\includegraphics[width=0.7\textwidth]{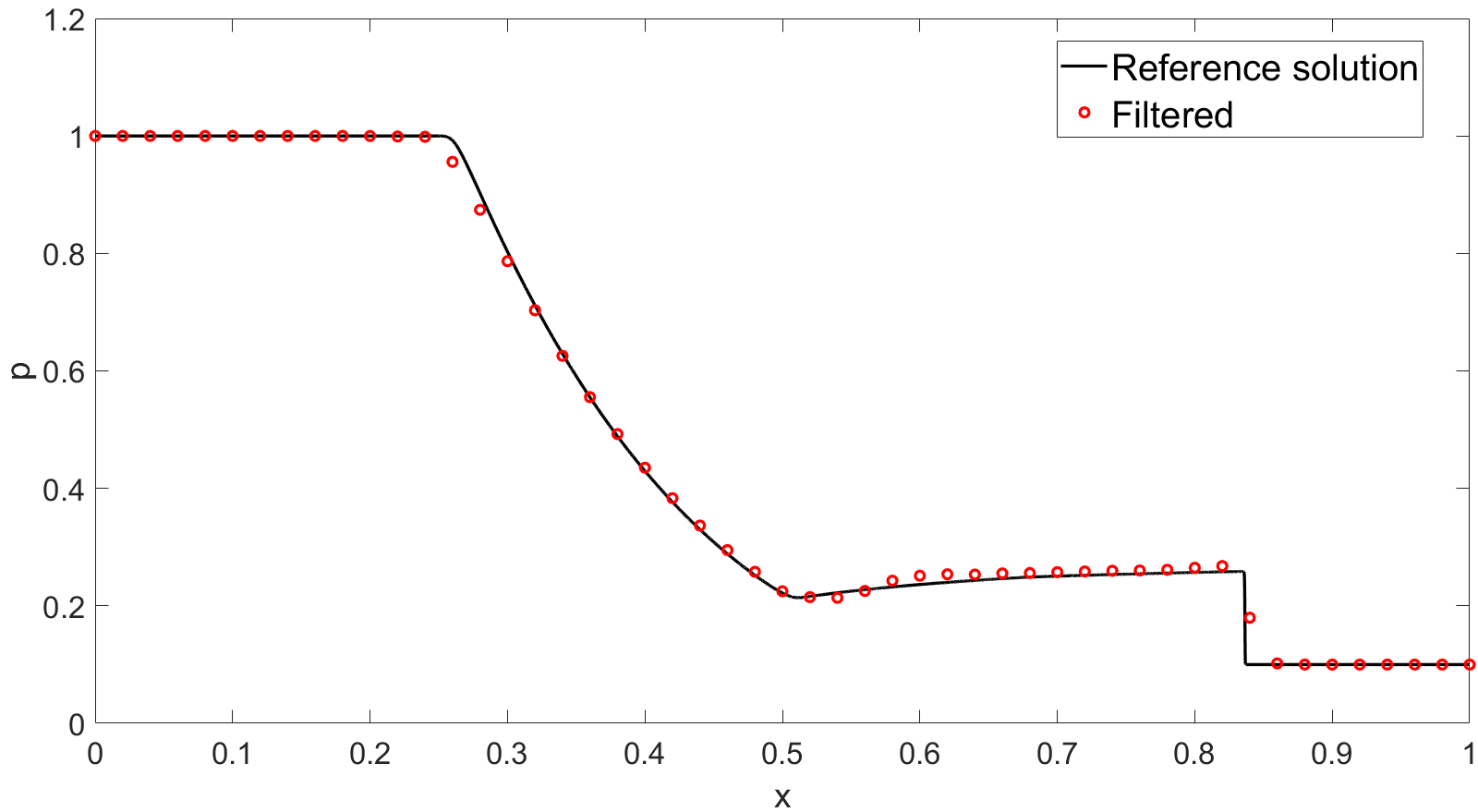} b)
	\end{subfigure}
	\caption{Computational results for 2D explosion problem at \(t = T_{f}\) with \(k = 1\), a) density, b) pressure. The black line reports the reference solution computed solving the \(1D\) problem in the radial direction, while the red dots represent the results of the simulation with Froese and Oberman's filter function using \(\beta_{\rho} = 1\), \(\beta_{\rho\mathbf{u}} = 1\) and \(\beta_{\rho E} = 1\).}
	\label{fig:2d_explosion_beta1_beta1_beta1}
\end{figure}

In this section, we consider the two-dimensional explosion problem discussed in \cite{dumbser:2014, zanotti:2015}. This test is quite relevant since it involves the propagation of waves that are not aligned with the mesh and therefore it can be used to check the ability of the proposed method to preserve physical symmetries of the problem as well as to validate it in multiple space dimensions. The computational domain is \(\Omega = \left(-1, 1\right)^2\), the final time is \(T_{f} = \SI{0.2}{\second}\) and the initial condition is the following:
\begin{equation}
\left(\rho_{0}, u_{0}, v_{0}, p_{0}\right) = \begin{cases}
\left(1, 0, 0, 1\right) \qquad &\text{if } r \le R \\
\left(0.125, 0, 0, 1\right) \qquad &\text{if } r > R,
\end{cases}
\end{equation} 
with \(R = 0.5\) denoting the radius of initial discontinuity and \(r = \sqrt{x^2 + y^2}\) representing the radial distance. As explained in \cite{toro:2009}, in \(2D\) we have cylindrical symmetry and a reference solution can be computed solving a one dimensional problem in the radial direction with suitable geometric source terms. Figure \ref{fig:2d_explosion_beta1_beta1_beta1} shows the results obtained using \(N_{el} = 200\) elements along each direction, \(k = 1\) and \(\beta_{\rho} = \beta_{\rho\mathbf{u}} = \beta_{\rho E} = 1\). One can easily notice that the discontinuities are well reproduced, even using only first order degree polynomial for the high order method, and their position is well captured with only slight undershoots and overshoots in correspondence of the rarefaction wave.    

The same test has been repeated increasing both the spatial resolution with \(N_{el} = 400\) and the high order polynomial degree with \(k = 2\). Figure \ref{fig:2d_explosion_beta1,7_beta1,7_beta1,7} reports the results obtained using \(\beta_{\rho} = \beta_{\rho\mathbf{u}} = \beta_{\rho E} = 1.7\) and an excellent agreement with the reference solution is achieved. Analogous results have been obtained in \cite{zanotti:2015}, where however polynomials of degree 9 were employed, and, for a 3D version of the problem, in \cite{loubere:2014}, with polynomials of degree 3.

\begin{figure}[pos = H]
	\begin{subfigure}[]{\textwidth}
		\centering
		\includegraphics[width=0.9\textwidth]{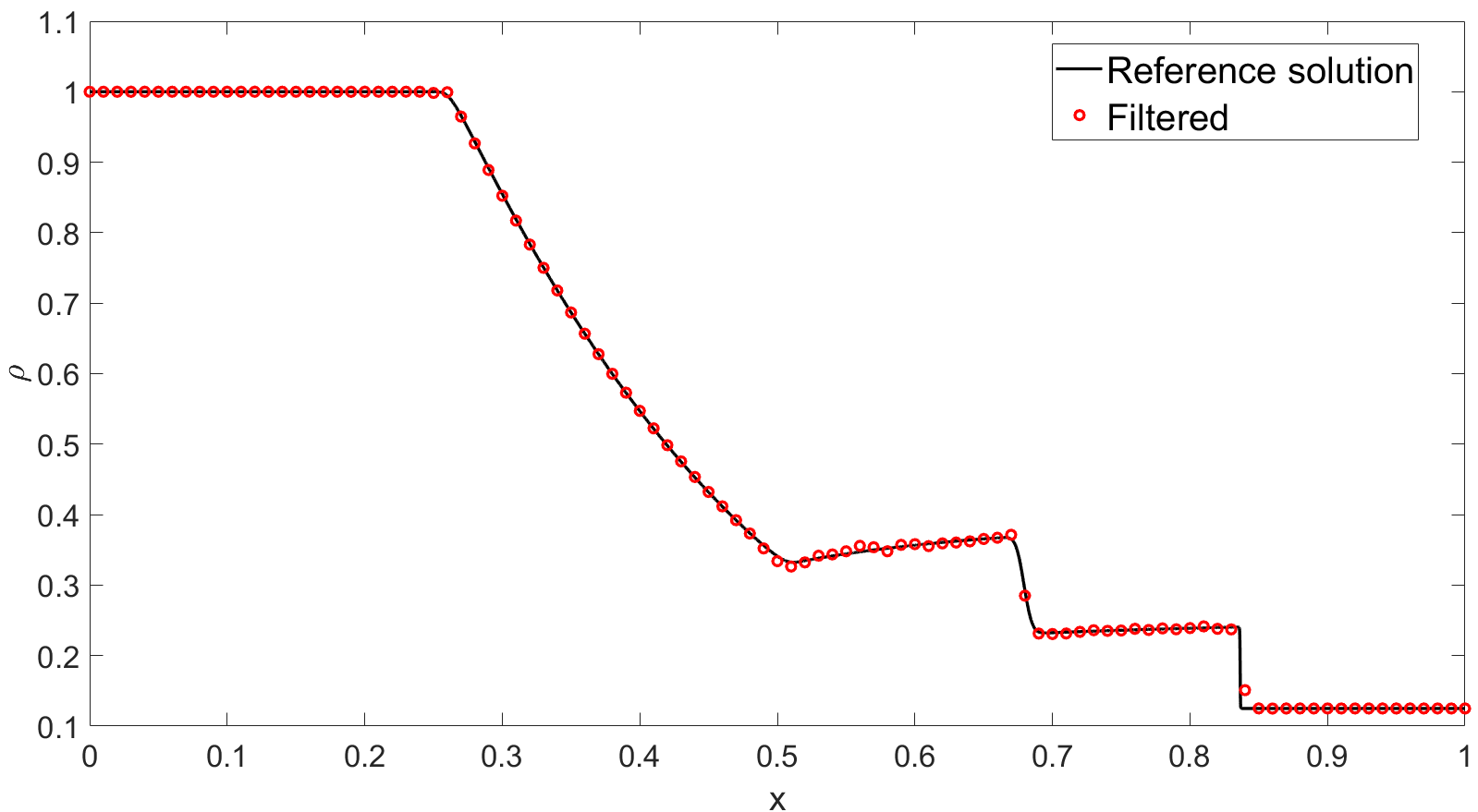} a) 
	\end{subfigure}
	\begin{subfigure}[]{\textwidth}
		\centering
		\includegraphics[width=0.9\textwidth]{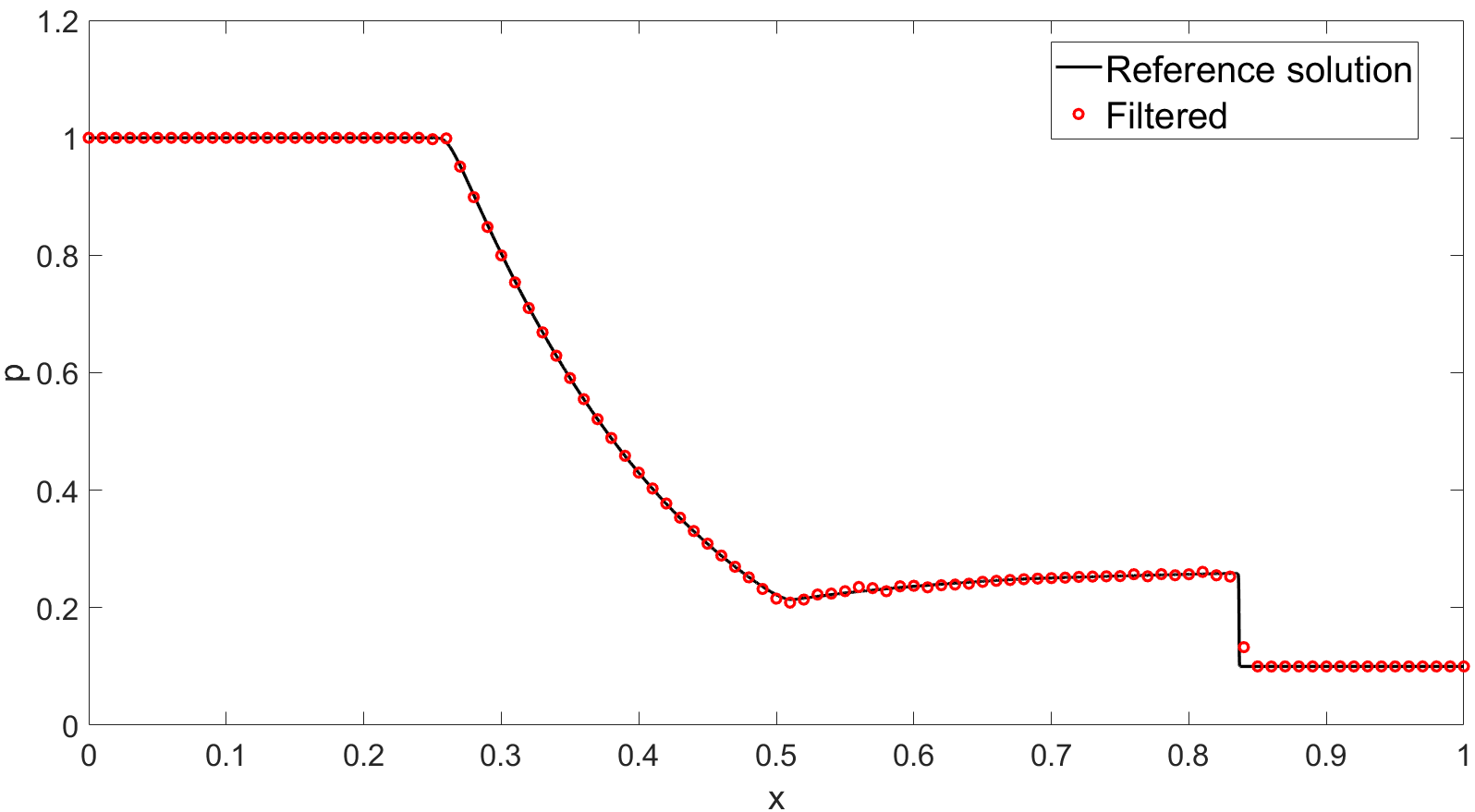} b)
	\end{subfigure}
	\caption{Computational results for 2D explosion problem at \(t = T_{f}\) with \(k = 2\), a) density, b) pressure. The black line reports the reference solution computed solving the \(1D\) problem in the radial direction, while the red dots represent the results of the simulation with Froese and Oberman's filter function using \(\beta_{\rho} = 1.7\), \(\beta_{\rho\mathbf{u}} = 1.7\) and \(\beta_{\rho E} = 1.7\).}
	\label{fig:2d_explosion_beta1,7_beta1,7_beta1,7}
\end{figure}

\begin{figure}[pos = H]
	\centering
	\includegraphics[width=0.9\textwidth]{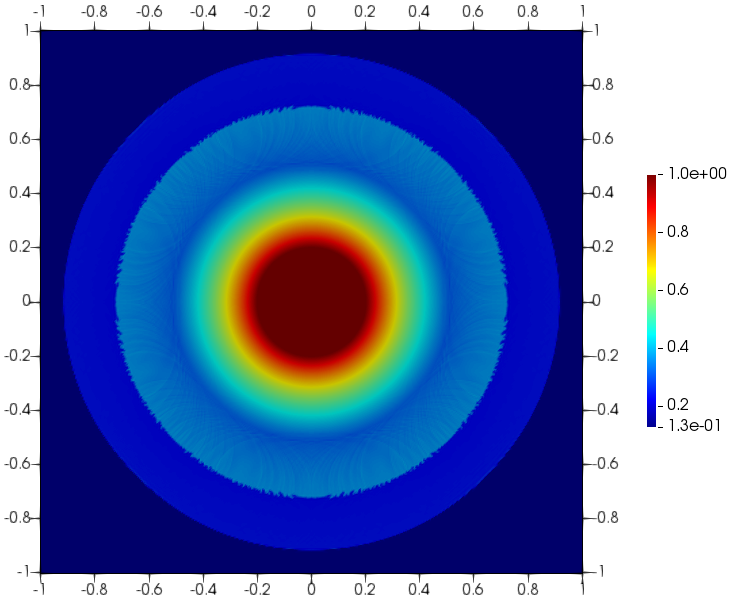} 
	\caption{Computational results for 2D explosion problem at \(t = T_{f}\) with \(k = 2\), contour plot of the density.}
	\label{fig:2d_explosion_beta1,7_beta1,7_beta1,7_2D}
\end{figure}

Finally, we have employed the \(h\)-adaptive version of the method, starting from a coarse mesh with \(N_{el} = 200\) elements along each direction and allowing up to three local refinements which would correspond to a uniform grid with \(N_{el} = 1600\). The employed local indicator is based on the gradient of the density; more specifically we define for each element \(K\) 
\begin{equation}\label{eq:adaptive_criterion_rho}
\eta_{K} = \max_{i \in \mathcal{N}_{K}} \left|\nabla \rho\right|_{i}.
\end{equation}
Figure \ref{fig:2d_explosion_beta1,7_beta1,7_beta1,7_adaptive_grid} shows the final grid obtained at \(t = T_{f}\) composed by 63136 elements and one can easily notice that more resolution is added in correspondence of the discontinuities.

\begin{figure}[pos = H]
	\begin{subfigure}[]{\textwidth}
		\centering
		\includegraphics[width=0.9\textwidth]{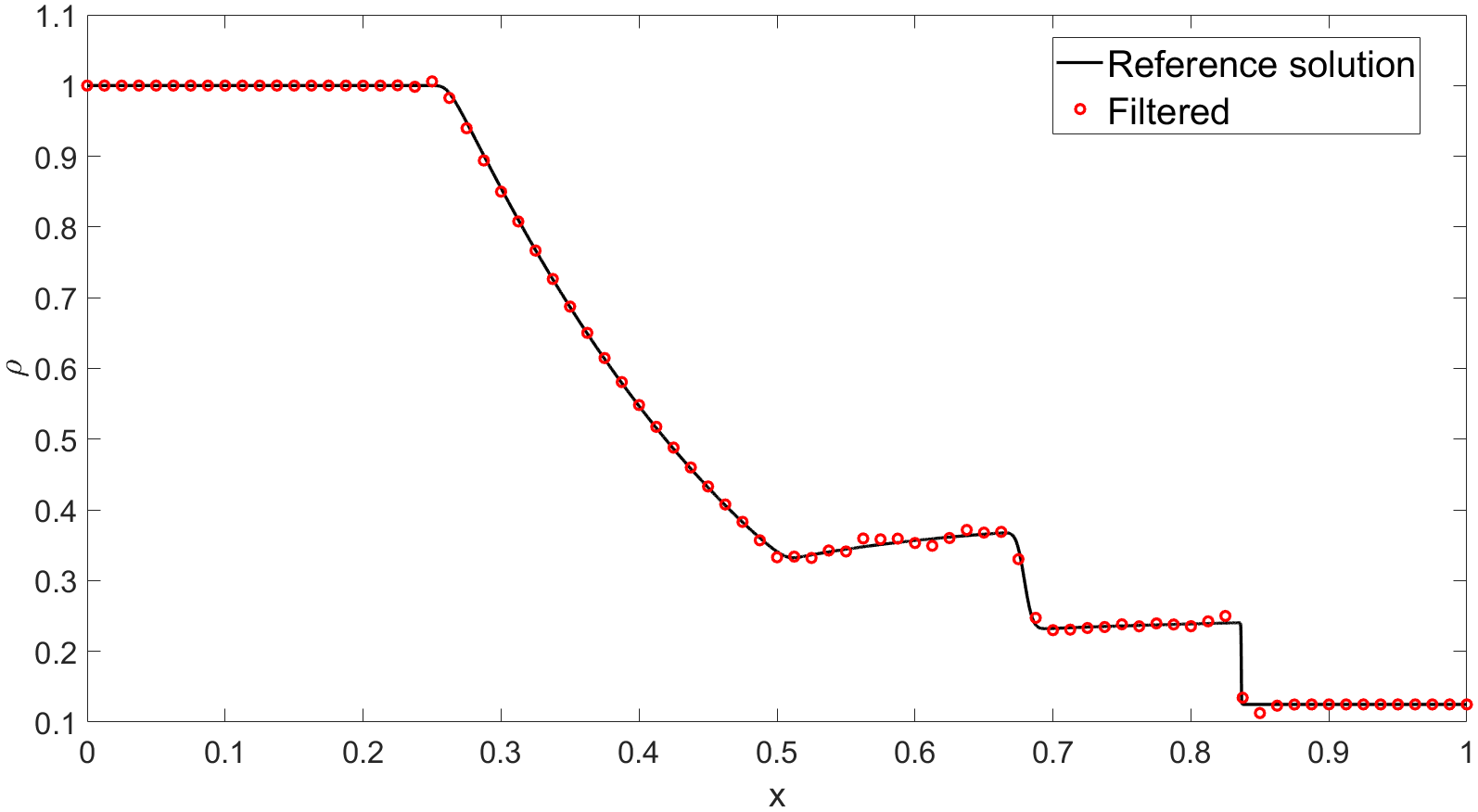} a) 
	\end{subfigure}
	\begin{subfigure}[]{\textwidth}
		\centering
		\includegraphics[width=0.9\textwidth]{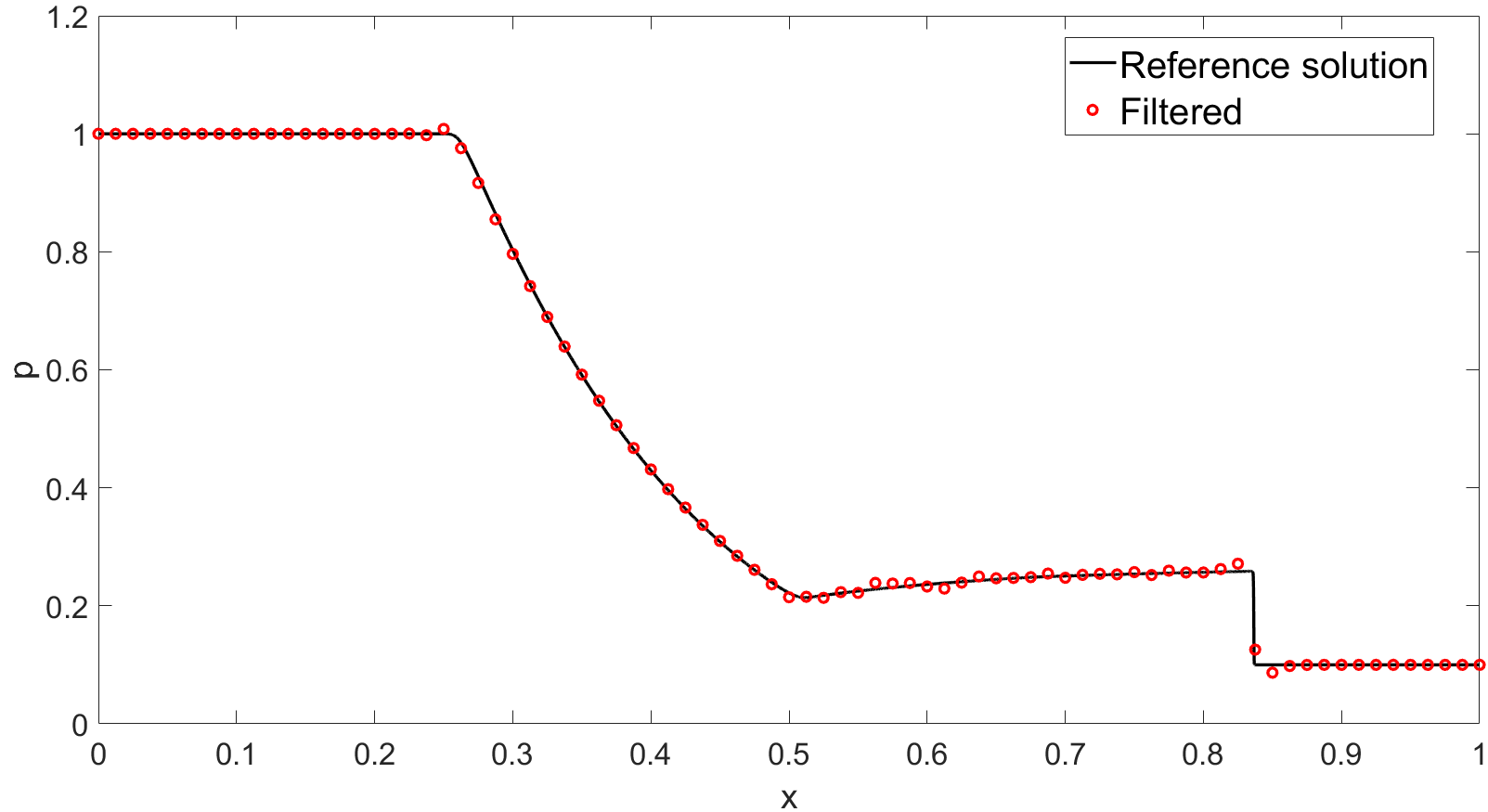} b)
	\end{subfigure}
	\caption{Computational results for 2D explosion problem at \(t = T_{f}\) with adaptive grid and \(k = 2\), a) density, b) pressure. The black line reports the reference solution computed solving the \(1D\) problem in the radial direction, while the red dots represent the results of the simulation with Froese and Oberman's filter function using \(\beta_{\rho} = 1.7\), \(\beta_{\rho\mathbf{u}} = 1.7\) and \(\beta_{\rho E} = 1.7\).}
	\label{fig:2d_explosion_beta1,7_beta1,7_beta1,7_adaptive}
\end{figure}

\begin{figure}[pos = H]
	\centering
	\includegraphics[width=0.9\textwidth]{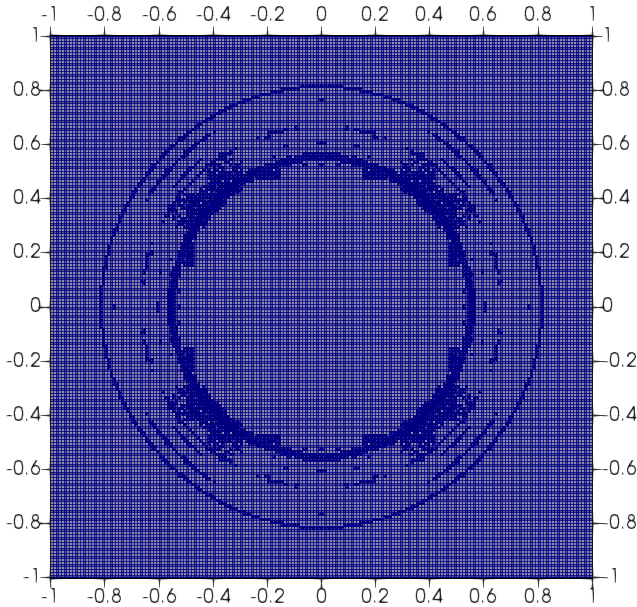}
	\caption{2D explosion problem, adaptive grid at \(t = T_{f}\) with \(k = 2\) and Froese and Oberman's filter function using \(\beta_{\rho} = 1.7\), \(\beta_{\rho\mathbf{u}} = 1.7\) and \(\beta_{\rho E} = 1.7\).}
	\label{fig:2d_explosion_beta1,7_beta1,7_beta1,7_adaptive_grid}
\end{figure}

\subsection{2D Riemann problem}
\label{ssec:2D_Riemann}

\begin{figure}[pos = H]
	\centering
	\includegraphics[width=0.8\textwidth]{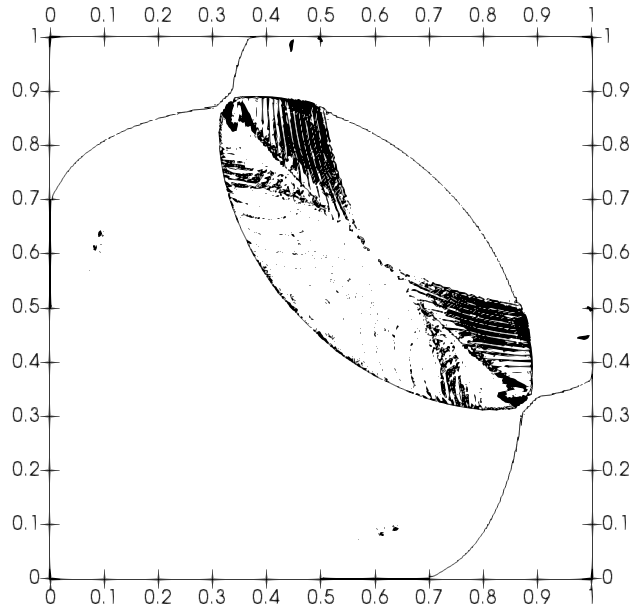}
	\caption{2D Riemann problem, isolines of the density at \(t = T_{f}\) with \(k = 2\) and Froese and Oberman's filter function using \(\beta_{\rho} = 0.25\), \(\beta_{\rho\mathbf{u}} = 0.25\) and \(\beta_{\rho E} = 0.25\).}
	\label{fig:2d_Riemann_rho}
\end{figure}

In this section, we consider a 2D Riemann problem corresponding to the Configuration 4 proposed in \cite{kurganov:2002}, which we summarize here for the convenience of the reader. The computational domain is \(\Omega = \left(0, 1\right)^{2}\) and the initial conditions are given by
\begin{equation}
\left(\rho_{0}, u_{0}, v_{0}, p_{0}\right) = \begin{cases}
\left(1.1, 0, 0, 1\right) \qquad &\text{if } x > 0.5 \text{ and } y > 0.5 \\
\left(0.5065, 0.8939, 0, 0.35\right) \qquad &\text{if } x < 0.5 \text{ and } y > 0.5 \\
\left(1.1, 0.8939, 0.8939, 1.1\right) \qquad &\text{if } x < 0.5 \text{ and } y < 0.5 \\
\left(0.5065, 0, 0.8939, 0.35\right) \qquad &\text{if } x > 0.5 \text{ and } y < 0.5.
\end{cases}
\end{equation}
The final time is \(T_{f} = \SI{0.25}{\second}\). In view of the particularly challenging conditions, we employ adaptive mesh refinement with the indicator described in \eqref{eq:adaptive_criterion_rho} in order to enhance the resolution along strong discontinuities. The initial mesh is composed by 200 elements along each direction and we allow up to two local refinements. We consider as high order polynomial degree \(k = 2\). Figure \ref{fig:2d_Riemann_rho} shows the results obtained for the density using \(\beta_{\rho} = \beta_{\rho\mathbf{u}} = \beta_{\rho E} = 0.25\). The filter tends to add more dissipation than needed, but this is necessary in order to avoid large undershoots and overshoots and more in general oscillations which completely corrupt the unfiltered solution. While not optimal, the results highlight the robustness of the proposed approach and show that the primary goal of the filter, namely avoid or at least reduce the oscillations, is achieved. Moreover, as pointed out in \cite{zanotti:2015}, the effects of Kelvin-Helmholtz instability with several small-scale features emerge at high resolution along the diagonal of the cocoon structure and this confirms that the test is particularly challenging. 

\section{Conclusions and future developments}
\label{sec:conclu}

In this work, a filtering technique for obtaining a monotonic Discontinuous Galerkin discretization of hyperbolic equations has been presented. The scheme is inspired by the approach originally proposed in \cite{bokanowski:2016} and it is based on a filter function that keeps the high order solution if it is regular and switches to a monotone low order approximation otherwise, according to the value of one or more parameters. Its potential has been demonstrated in a number of classical benchmarks for linear advection and Euler equations. 

In future work, an obvious and necessary development concerns the tuning of the parameter(s) \(\beta\). The goal is to automatically choose suitable values depending on the employed time and space steps as well as the polynomial degree used by the higher order discretization. Moreover, we plan to investigate the behaviour of the proposed method in case stiff source terms and/or non-conservative terms are present, as for example in the Baer-Nunziato model of compressible multiphase flows \cite{baer:1986}. 

\section*{Acknowledgements}
The author would like to thank Luca Bonaventura for several useful discussions. The author also greatfully acknowledges the two anonymous reviewers, which have greatly helped in improving the quality of the paper.

\bibliographystyle{cas-model2-names}

\bibliography{DG_monotonization.bib}

\end{document}